\documentclass[journal]{IEEEtran}

\usepackage{url}
\usepackage{graphicx}
\usepackage{subfigure}
\usepackage{color}
\usepackage{epsfig}
\usepackage{amssymb}
\usepackage{latexsym}

\newcommand {\C} {{\rm I\kern-5.5pt C}}

\newcommand{\bP}[1]{{\mathbb{P}}\left[{#1}\right]}
\newcommand{\bE}[1]{{\mathbb{E}}\left[{#1}\right]}

\newcommand{\1}[1]{{\bf 1}\left[#1\right]}       

\newcommand{\fsquare}{\vrule height6pt width7pt depth1pt}   
\newcommand{\myproof}{{\hfill \\ \bf Proof. \ }}           
\newcommand{\myendpf}{\hfill\fsquare \\[0.1in]}             

\newtheorem{theorem}{Theorem}[section]

\newtheorem{lemma}[theorem]{Lemma}
\newtheorem{proposition}[theorem]{Proposition}

\begin{document}

\title{Zero-one laws for connectivity in \\random key graphs
\thanks{Manuscript received November 16, 2010; revised August 22, 2011.
This work was supported by NSF Grant CCF-07290. The
material in this paper
was presented in part at
the 2008 IEEE International Symposium
on Information Theory (ISIT 2008), Toronto (Canada), June 2008,
and at the 2009
IEEE International Symposium on Information Theory (ISIT 2009),
Seoul (S. Korea), June 2009.}\thanks{O. Ya\u{g}an 
was with the Department of
Electrical and Computer Engineering, and the Institute for Systems
Research, University of Maryland, College Park, MD 20742 USA. He
is now with CyLab, Carnegie Mellon University, 
Pittsburgh, PA 15213 USA (e-mail: osmanyagan@gmail.com).} 
\thanks{A. M. Makowski is with the
Department of Electrical and Computer Engineering, and the
Institute for Systems Research, University of Maryland, College
Park, MD 20742 USA (e-mail: armand@isr.umd.edu).}
\thanks{Copyright (c) 2011 IEEE. Personal use of this material is
permitted.  However, permission to use this material for any other
purposes must be obtained from the IEEE by sending a request to
pubs-permissions@ieee.org.}}

\author{
Osman Ya\u{g}an and Armand M. Makowski, \it{Fellow, IEEE}
}

\date{\today}
\maketitle

\begin{abstract}
\normalsize The random key graph is a random graph naturally
associated with the random key predistribution scheme introduced by
Eschenauer and Gligor in the context of wireless sensor networks. For this class
of random graphs we establish a new version of a conjectured
zero-one law for graph connectivity as the number of nodes becomes
unboundedly large. The results reported here complement and
strengthen recent work on this conjecture by Blackburn and Gerke.
In particular, the results are given under conditions which are
more realistic for applications to wireless sensor networks.
\end{abstract}


\section{Introduction}
\label{sec:Introduction}

\subsection{Background}
\label{subsec:IntroductionBackground}

Random key graphs, also known as uniform random intersection
graphs, are random graphs that belong to the class of random
intersection graphs \cite{SingerThesis}. They have appeared
recently in application areas as diverse as clustering analysis
\cite{GodehardtJaworski, GodehardtJaworskiRybarczyk},
collaborative filtering in recommender systems \cite{Marbach2008}
and random key predistribution for wireless sensor networks (WSNs)
\cite{DiPietroManciniMeiPanconesiRadhakrishnan2006,
DiPietroManciniMeiPanconesiRadhakrishnan2008,EschenauerGligor}.

For the sake of concreteness, we introduce this class of random
graphs in this last context (hence the terminology). A WSN is a
collection of spatially distributed sensors with limited
capabilities for computations and wireless communications. It is
envisioned that such networks will be used in applications such as
battlefield surveillance, environment monitoring and traffic
control, to name a few. In many settings, both military and
civilian, network security will be a basic requirement for
successful operations. However, traditional key exchange and distribution
protocols are based on trusting third parties, and turn out to be
inadequate for large-scale wireless sensor networks, e.g., see
\cite{EschenauerGligor,PerrigStankovicWagner,SunHe,WangAtteburyRamamurthy}
for discussions of some of the challenges. To address some of the
difficulties Eschenauer and Gligor \cite{EschenauerGligor} have
recently proposed the following random key predistribution scheme:

Before deployment, each sensor in a WSN is independently assigned
$K$ distinct cryptographic keys which are selected at random from
a pool of $P$ keys (with $K < P$). These $K$ keys constitute the
key ring of the node and are inserted into its memory. Two sensor
nodes can then establish a secure link between them if they are
within transmission range of each other {\em and} if their key
rings have at least one key in common; see \cite{EschenauerGligor}
for implementation details. A situation of particular interest is
that of {\em full visibility} whereby nodes are all within
communication range of each other. In that case a secure link can be
established between two nodes if their key rings have at least one
key in common. The  resulting notion of adjacency defines
the {\em random key} graph $\mathbb{K}(n;(K,P))$ on the vertex set 
$\{ 1, \ldots , n \}$ where $n$ is the number of sensor nodes; see
Section \ref{sec:RandomKeyGraph} for precise definitions.

A basic question concerning the scheme of Eschenauer and Gligor is
its ability to achieve {\em secure connectivity} amongst
participating nodes in the sense that a {\em secure path} exists
between any pair of nodes. Therefore, under full visibility it is
natural to seek conditions on $n$, $K$ and $P$ under which the
random key graph $\mathbb{K}(n;(K,P))$ constitutes a connected
graph with high probability -- The availability of such conditions
would provide an encouraging indication of the feasibility of
using this distribution scheme for WSNs. As discussed in Section
\ref{sec:BasisConjecture}, this search has lead to {\em
conjecturing} the following zero-one law for graph connectivity in
random key graphs: If the parameters $K$ and $P$ are scaled with
$n$ according to
\begin{equation}
\frac{K_n^2}{P_n} = \frac{ \log n + \alpha _n }{n}, \quad n=1,2,
\ldots \label{eq:CONJECTUREScalingK+P}
\end{equation}
for some sequence $\alpha: \mathbb{N}_0 \rightarrow \mathbb{R}$,
then it has been conjectured that
\begin{eqnarray}
\lefteqn{ \lim_{n \rightarrow \infty } \bP{ \mathbb{K}(n;
(K_n,P_n) ) ~{\rm is~connected} } } & &
\nonumber \\
&=& \left \{
\begin{array}{ll}
0 & \mbox{if~ $\lim_{ n\rightarrow \infty }\alpha_n = - \infty $} \\
  &                 \\
1 & \mbox{if~ $\lim_{ n\rightarrow \infty }\alpha_n = + \infty $.}
\end{array}
\right . \label{eq:CONJECTUREZeroOneLawConnectivity}
\end{eqnarray}
This conjecture appeared independently in
\cite{BlackburnGerke,YaganMakowskiISIT2008}. The zero-one law
(\ref{eq:CONJECTUREScalingK+P})-(\ref{eq:CONJECTUREZeroOneLawConnectivity})
mimics a similar one for Erd\H{o}s-R\'enyi graphs \cite{Bollobas},
and can be motivated from it by asymptotically matching the link
assignment probabilities in these two classes of random graphs.

\subsection{Related work}
\label{subsec:RelatedWork}

Recent results concerning the conjectured zero-one law
(\ref{eq:CONJECTUREScalingK+P})-(\ref{eq:CONJECTUREZeroOneLawConnectivity})
are now surveyed: Di Pietro et al. have shown \cite[Thm.
4.6]{DiPietroManciniMeiPanconesiRadhakrishnan2008} that for large
$n$, the random key graph will be connected with very high
probability if $P_n$ and $K_n$ are selected such that
\[
K_n \geq 5, \ P_n \geq n \quad \mbox{and} \quad \frac{K_n^2}{P_n}
\sim c~ \frac{\log n}{n} \label{eq:Selection}
\]
as soon as $c \geq 16$.\footnote{In the conference version of this
work \cite[Thm. 4.6]{DiPietroManciniMeiPanconesiRadhakrishnan2006}
the result is claimed to hold for $c > 8$.} They also observe that
for large $n$, the random key graph will be disconnected with very
high probability if the scaling satisfies
\[
\frac{K^2_n}{P_n} = o \left ( \frac{\log n}{n} \right ).
\]

The zero-law in (\ref{eq:CONJECTUREZeroOneLawConnectivity}) has
recently been established independently by Godehardt and Jaworski
\cite{GodehardtJaworski}, Blackburn and Gerke
\cite{BlackburnGerke}, and Ya\u{g}an and Makowski
\cite{YaganMakowskiISIT2008}. In all these papers, it was shown that
\[
\lim_{n \to \infty} \bP{\mathbb{K}(n; (K_n,P_n) ) ~{\rm
contains~no~isolated~nodes} } = 0
\]
whenever $\lim_{n \to \infty} \alpha_n = -\infty$ in
(\ref{eq:CONJECTUREScalingK+P}), a result which clearly implies
the conjectured zero-law.

Blackburn and Gerke \cite{BlackburnGerke} also succeeded in
generalizing the one-law result by Di Pietro et al. in a number of
directions: Under the additional conditions
\begin{equation}
K_n \geq 2 \quad \mbox{and} \quad  P_n \geq n, \quad n=1,2,
\ldots, \label{eq:OneBGMainCondition}
\end{equation}
they showed \cite[Thm. 5]{BlackburnGerke} that
\begin{equation}
\lim_{n \rightarrow \infty} \bP{\mathbb{K}(n; (K_n,P_n) ) ~{\rm
is~connected} } = 1 \label{eq:OneBGMainResult1}
\end{equation}
if
\begin{equation}
\liminf_{n \rightarrow \infty} \frac{K^2_n}{P_n} \frac{n}{\log n}
> 1. \label{eq:OneBGMainResult2}
\end{equation}
This result is weaker than the one-law in the conjecture
(\ref{eq:CONJECTUREScalingK+P})-(\ref{eq:CONJECTUREZeroOneLawConnectivity}).
However, in the process of establishing
(\ref{eq:OneBGMainResult1})-(\ref{eq:OneBGMainResult2}), they also
show \cite[Thm. 3]{BlackburnGerke} that the conjecture does hold
in the special case $K_n = 2$ for all $n=1,2, \ldots $ {\em
without} any constraint on the size of the key pools, say $P_n
\leq n $ or $n \leq P_n$. Specifically, the one-law in
(\ref{eq:CONJECTUREZeroOneLawConnectivity}) is shown to hold
whenever the scaling is done according to
\[
K_n =2 , \quad \frac{4}{P_n} = \frac{ \log n + \alpha_n}{n}, \quad
n= 1,2,\ldots
\]
as soon as $\lim_{n \rightarrow \infty} \alpha_n = \infty $. As
pointed out by these authors, it is now easy to conclude that the
one-law in (\ref{eq:CONJECTUREZeroOneLawConnectivity}) holds
whenever $2 \leq K_n \leq P_n$ and $P_n = o \left ( \frac{n}{\log
n} \right )$; this corresponds to a constraint $P_n \ll n $.

\subsection{Contributions}
\label{subsec:Contributions}

In this paper, we complement existing results concerning the
conjecture
(\ref{eq:CONJECTUREScalingK+P})-(\ref{eq:CONJECTUREZeroOneLawConnectivity})
in several ways: We establish (Theorem
\ref{thm:ZeroOneLaw+Connectivity}) the one-law in
(\ref{eq:CONJECTUREZeroOneLawConnectivity}) under the conditions
$K_n \geq 2$ and $P_n=\Omega(n)$, i.e., $P_n \geq \sigma n$ for
some $\sigma>0$.  Since the zero-law in
(\ref{eq:CONJECTUREZeroOneLawConnectivity}) has already been
established \cite{BlackburnGerke, GodehardtJaworski,
YaganMakowskiISIT2008}, the validity of
(\ref{eq:CONJECTUREScalingK+P})-(\ref{eq:CONJECTUREZeroOneLawConnectivity})
thus follows whenever $P_n=\Omega(n)$ and $K_n \geq 2$.

This result already improves on the one-law
(\ref{eq:OneBGMainResult1})-(\ref{eq:OneBGMainResult2}) obtained
by Blackburn and Gerke \cite{BlackburnGerke} under the condition
(\ref{eq:OneBGMainCondition}). Moreover, as discussed earlier,
these authors have established the conjectured one-law in
(\ref{eq:CONJECTUREZeroOneLawConnectivity}) under conditions
very different from the ones used here, i..e., either
$K_n=2$ or $K_n \geq 2$ with $P_n = o \left ( \frac{n}{\log n}
\right)$. In practical WSN scenarios it is expected that the size
of the key pool will be much larger than the number of
participating nodes
\cite{DiPietroManciniMeiPanconesiRadhakrishnan2008,
EschenauerGligor} and that key rings will contain more than two
keys. In this context, our results concerning the full conjecture
(\ref{eq:CONJECTUREScalingK+P})-(\ref{eq:CONJECTUREZeroOneLawConnectivity})
are therefore given under more realistic conditions than earlier
work.

The proof of the main result is lengthy and technically involved.
However, in a parallel development, we have also shown in
\cite{YaganMakowskiCISS2010} that when $P_n = O(n^\delta)$ with $0
< \delta < \frac{1}{2}$, the so-called small key pool case,
elementary arguments can be used to establish a one-law for
connectivity. This is an easy byproduct of the observation that
connectivity is achieved in the random key graph whenever {\em
all} possible key rings have been distributed to the participating
nodes.

The results established in this paper were first announced in the
conference paper \cite{YaganMakowskiSeoul2009A} with an outline of
the proofs; the full details were provided in an early draft
\cite{YaganMakowskiISRTechReport2009} posted in January 2009.
However, after completing this work, we learned of the independent
work of Rybarczyk \cite{Rybarczyk2009} concerning the conjecture
(\ref{eq:CONJECTUREScalingK+P})-(\ref{eq:CONJECTUREZeroOneLawConnectivity})
without any condition on the size of the key pool. Reference
\cite{Rybarczyk2009} deals mainly with the diameter and phase
transition threshold of random key graphs, and uses branching
process arguments similar to the ones given in \cite{ChungLu}. The
intermediary results, the so-called branching process lemmas, pave
the way to a proof of the conjecture
(\ref{eq:CONJECTUREScalingK+P})-(\ref{eq:CONJECTUREZeroOneLawConnectivity})
by an approach very different from the one used here.

\subsection{The structure of the paper}

The paper is organized as follows: The class of random key graphs
is formally introduced in Section \ref{sec:RandomKeyGraph}. A
basis for the conjectured zero-one law is discussed in Section
\ref{sec:BasisConjecture}, and the main result of the paper,
summarized as Theorem \ref{thm:ZeroOneLaw+Connectivity}, is
presented in Section \ref{sec:MainResult}. A roadmap to the proof
of Theorem \ref{thm:ZeroOneLaw+Connectivity} is given in Section
\ref{sec:RoadMap}. The approach is similar to the one used for
proving the one-law for graph connectivity in Erd\H{o}s-R\'enyi
graphs \cite[p. 164]{Bollobas},
\cite[Section 3.4, p. 40]{DraiefMassoulieLMS} ,
\cite[p. 304]{SpencerSaintFlour1991}; see
(\ref{eq:ER+ScalingStrong})-(\ref{eq:ERZeroOneLawConnectivity}).
Here as well, we focus on the probability that the random key
graph is not connected and yet has no isolated nodes. We then seek
to show that this probability becomes vanishingly small as $n$
grows large under the appropriate scaling. As in the classical
case this is achieved through a combination of judicious bounding
arguments, the starting point being the well-known bound
(\ref{eq:BasicIdea+UnionBound2}) on the probability of interest.
However, in order for these arguments to successfully go through,
we found it necessary to restrict attention to a subclass of
structured scalings (referred throughout as strongly admissible
scalings). In Section \ref{sec:ReductionStep} a reduction argument
shows that we need only establish the desired one-law for such
strongly admissible scalings. The explanation of the right
handside of (\ref{eq:CONJECTUREScalingK+P}) as a proxy for link
assignment in the limiting regime is revealed through a useful
equivalence developed in Section \ref{sec:Usefulequivalence}.

With these technical prerequisites in place, the needed bounding
arguments are then developed in Section \ref{sec:BasicUnionBound},
Section \ref{sec:BoundingProbabilities} and Section
\ref{sec:Tail+ImprovedBounds}, and the final steps of the proof of
Theorem \ref{thm:ZeroOneLaw+Connectivity} are outlined in Section
\ref{sec:Outline}. The final sections of the paper, namely Section
\ref{sec:ToShowSecondPiece0} through Section
\ref{sec:ProofPropositionRef{prop:D}}, are devoted to the various
technical steps needed to complete the arguments outlined in
Section \ref{sec:Outline}.

\subsection{Notation and conventions}
\label{subsec:Notation}

A word on the notation and conventions in use: All limiting
statements, including asymptotic equivalences, are understood with
$n$ going to infinity. The random variables (rvs) under
consideration are all defined on the same probability triple
$(\Omega, {\cal F}, \mathbb{P})$. Probabilistic statements are
made with respect to this probability measure $\mathbb{P}$, and we
denote the corresponding expectation operator by $\mathbb{E}$. The
indicator function of an event $E$ is denoted by $\1{E}$. For any
discrete set $S$ we write $|S|$ for its cardinality.

\section{Random key graphs}
\label{sec:RandomKeyGraph}

Random key graphs are parametrized by the number $n$ of nodes, the
size $P$ of the key pool and the size $K$ of each key ring with $K
\leq P$. To lighten the notation we often group the integers $P$
and $K$ into the ordered pair $\theta \equiv (K,P)$.

Nodes are labelled $i, \ldots , n$ while keys are labelled $1,
\ldots , P$. For each node $i=1, \ldots , n$, let $K_i (\theta)$
denote the random set of $K$ distinct keys assigned to node $i$.
We can think of $K_i(\theta)$ as an $\mathcal{P}_{K} $-valued rv
where $\mathcal{P}_{K} $ denotes the collection of all subsets of
$\{ 1, \ldots , P \}$ which contain exactly $K$ elements --
Obviously, we have $|\mathcal{P}_{K} | =  {P \choose K}$. The rvs
$K_1(\theta), \ldots , K_n(\theta)$ are assumed to be {\em i.i.d.}
rvs, each of which is {\em uniformly} distributed over
$\mathcal{P}_{K}$ with
\[
\bP{ K_i(\theta) = S } = {P \choose K} ^{-1}, \quad S \in
\mathcal{P}_{K} \label{eq:KeyDistrbution1}
\]
for all $i=1, \ldots , n$. This corresponds to selecting keys
randomly and {\em without} replacement from the key pool.

Distinct nodes $i,j=1, \ldots , n$ are said to be adjacent if they
share at least one key in their key rings, namely
\[
K_i (\theta) \cap K_j (\theta) \not = \emptyset,
\label{eq:KeyGraphConnectivity}
\]
in which case an undirected link is assigned between nodes $i$ and
$j$. The resulting random graph defines the {\em random key graph}
on the vertex set $\{ 1, \ldots , n\}$, hereafter denoted by
$\mathbb{K}(n; \theta )$. For distinct $i,j =1, \ldots , n$, it is
a simple matter to check that
\[
\bP{ K_i (\theta) \cap K_j (\theta) = \emptyset } = q (\theta)
\]
with
\begin{equation}
q (\theta) = \left \{
\begin{array}{ll}
0 & \mbox{if~ $P <2K$} \\
  &                 \\
\frac{{P-K \choose K}}{{P \choose K}} & \mbox{if~ $2K \leq P$,}
\end{array}
\right . \label{eq:q_theta}
\end{equation}
whence the probability of edge occurrence between any two nodes is
equal to $1-q(\theta)$. The expression (\ref{eq:q_theta}) and
others given later are simple consequences of the often used fact
that
\begin{equation}
\bP{ S \cap K_i(\theta) = \emptyset } = \left \{
\begin{array}{ll}
0 & \mbox{if~ $|S| > P-K$} \\
  &                 \\
\frac{{P- |S| \choose K}}{{P \choose K}} & \mbox{if~ $|S| \leq
P-K$}
\end{array}
\right . \label{eq:Probab_key_ring_does_not_intersect_S}
\end{equation}
with $S$ a subset of $\{ 1, \ldots , P \}$. The case $P<2K$
corresponds to an edge existing between every pair of nodes, so
that $\mathbb{K}(n;\theta)$ coincides with the complete graph on
the vertex set $\{ 1, \ldots , n \}$. Also, we always have $0 \leq
q(\theta) < 1 $ with $q(\theta)> 0$ if and only if $2K \leq P$.

Random key graphs form a subclass in the family of {\em random
intersection} graphs. However, the model adopted here differs from
the random intersection graphs discussed by Singer-Cohen et al. in
\cite{KaronskiScheinermanSingerCohen,SingerThesis} where each node
is assigned a key ring, one key at a time according to a
Bernoulli-like mechanism (so that each key ring has a random size
and has positive probability of being empty). Both subclasses are
subsumed by the more general random intersection graph model
discussed by Godehardt et al. \cite{GodehardtJaworski,
GodehardtJaworskiRybarczyk}.

Throughout, with $n=2,3, \ldots $, and positive integers $K$ and
$P$ such that $K \leq P$, let $P(n;\theta)$ denote the probability
that the random key graph $\mathbb{K}(n; \theta)$ is connected,
namely
\[
P(n;\theta) := \bP{ \mathbb{K}(n; \theta) ~{\rm is~connected} },
\quad \theta = (K,P).
\]

\section{A basis for the conjecture}
\label{sec:BasisConjecture}

As indicated earlier, we wish to select $P$ and $K$ so that
$P(n;\theta)$ is as large (i.e., as close to one) as possible. We
outline below a possible approach which is inspired by the
discussion on this issue given by Eschenauer and Gligor in their
original work \cite{EschenauerGligor}; see also the discussion in
\cite{DiPietroManciniMeiPanconesiRadhakrishnan2006,
DiPietroManciniMeiPanconesiRadhakrishnan2008},

(i) Let $\mathbb{G}(n;p)$ denote the Erd\H{o}s-R\'enyi graph on
$n$ vertices with edge probability $p$ ($0< p \leq 1$)
\cite{Bollobas, DraiefMassoulieLMS, JansonLuczakRucinski}. Despite
strong similarities, the random graph $\mathbb{K}(n; \theta)$ is
{\em not} an Erd\H{o}s-R\'enyi graph $\mathbb{G}(n;p)$. This is so
because edge assignments are independent in $\mathbb{G}(n;p)$ but
can be correlated in $\mathbb{K}(n; \theta)$. Yet, setting aside
this (inconvenient) fact, we note that $\mathbb{K}(n;\theta)$ can
be matched naturally to an Erd\H{o}s-R\'enyi graph
$\mathbb{G}(n;p)$ with $p$ and $\theta$ related through
\begin{equation}
p = 1 - q(\theta). \label{eq:Proxy}
\end{equation}
This constraint ensures that link assignment probabilities in
$\mathbb{K}(n;\theta)$ and $\mathbb{G}(n;p)$ coincide. Moreover,
under (\ref{eq:Proxy}) it is easy to check that the degree of a
node in either random graph is a Binomial rv with the same
parameters, namely $n-1$ and $p=1-q(\theta)$!\footnote{For
Erd\H{o}s-R\'enyi graphs this result is well known, while for
random key graphs this characterization is a straightforward
consequence of (\ref{eq:Probab_key_ring_does_not_intersect_S}).}
Given that the degree distributions in a random graph are often
taken (perhaps mistakenly) as a good indicator of its connectivity
properties, it is tempting to conclude that the zero-one law for
graph connectivity in random key graphs can be inferred from the
analog result for Erd\H{o}s-R\'enyi graphs when matched through
the condition (\ref{eq:Proxy}).

(ii) To perform such a \lq\lq transfer," we first recall that in
Erd\H{o}s-R\'enyi graphs the property of graph connectivity is
known to exhibit the following zero-one law \cite{Bollobas}: If we
scale the edge assignment probability $p$ according to
\begin{equation}
p_n = \frac{ \log n + \alpha _n }{n}, \quad n=1,2, \ldots
\label{eq:ER+ScalingStrong}
\end{equation}
for some sequence $\alpha: \mathbb{N}_0 \rightarrow \mathbb{R}$,
then
\begin{eqnarray}
\lefteqn{ \lim_{n \rightarrow \infty } \bP{ \mathbb{G}(n; p_n )
~{\rm is~connected} } } & &
\nonumber \\
&=& \left \{
\begin{array}{ll}
0 & \mbox{if~ $\lim_{ n\rightarrow \infty }\alpha_n = - \infty $} \\
  &                 \\
1 & \mbox{if~ $\lim_{ n\rightarrow \infty }\alpha_n = + \infty $.}
\end{array}
\right . \label{eq:ERZeroOneLawConnectivity}
\end{eqnarray}

(iii) Under the matching condition (\ref{eq:Proxy}), these
classical results suggest scaling the parameters $K$ and $P$ with
$n$ according to
\begin{equation}
1 - \frac{ {P_n-K_n \choose K_n} }{ {P_n \choose K_n} } = \frac{
\log n + \alpha _n }{n}, \quad n=1,2, \ldots \label{eq:ScalingK+P}
\end{equation}
for some sequence $\alpha: \mathbb{N}_0 \rightarrow \mathbb{R}$.
In view of (\ref{eq:ERZeroOneLawConnectivity}) it is then not too
unreasonable to expect that the zero-one law
\begin{eqnarray}
\lim_{n \rightarrow \infty } P(n;\theta_n) = \left \{
\begin{array}{ll}
0 & \mbox{if~ $\lim_{ n\rightarrow \infty }\alpha_n = - \infty $} \\
  &                 \\
1 & \mbox{if~ $\lim_{ n\rightarrow \infty }\alpha_n = + \infty $}
\end{array}
\right . \label{eq:Pre-ConjecturedZeroOneLawConnectivity}
\end{eqnarray}
should hold (possibly under some additional assumptions).

Of course, for this approach to be operationally useful, a good
approximation to the right handside of (\ref{eq:Proxy}) is needed.
Eschenauer and Gligor provided such an approximation with the help
of Stirling's formula. However, as already indicated by Di Pietro
et al. \cite{DiPietroManciniMeiPanconesiRadhakrishnan2006,
DiPietroManciniMeiPanconesiRadhakrishnan2008}, it is easy to check
that
\begin{equation}
1 - \frac{ {P-K \choose K} }{ {P \choose K } } \simeq
\frac{K^2}{P} \label{eq:EquivalenceNeeded}
\end{equation}
under natural assumptions; see Lemma
\ref{lem:AsymptoticEquivalence}. Thus, if instead of scaling the
parameters according to (\ref{eq:ScalingK+P}), we scale them
according to
\[
\frac{K_n^2}{P_n} = \frac{ \log n + \alpha _n }{n}, \quad n=1,2,
\ldots \label{eq:ScalingK+PAlternative}
\]
then it is natural to conjecture that the zero-one law
(\ref{eq:Pre-ConjecturedZeroOneLawConnectivity}) should still
hold.

While this transfer technique could in principle be applied to
other graph properties, it may not always yield the correct form
for the zero-one law; see  the papers
\cite{YaganMakowskiAllerton2009, YaganMakowskiTriangle} for
results concerning the existence of triangles in random key graphs.

\section{The main result}
\label{sec:MainResult}

Any pair of functions $P,K: \mathbb{N}_0 \rightarrow \mathbb{N}_0$
defines a {\em scaling} provided the natural conditions
\[
K_n \leq P_n , \quad n=1,2, \ldots \label{eq:AdmissibilityC}
\]
are satisfied. We can always associate with it a sequence 
$\alpha : \mathbb{N}_0 \rightarrow \mathbb{R}$ 
through the relation
\begin{equation}
\frac{K_n^2}{P_n} = \frac{ \log n + \alpha_n }{n}, \quad n=1, 2,
\ldots \label{eq:DeviationCondition}
\end{equation}
Just set
\[
\alpha_n := n \frac{K_n^2}{P_n} - \log n, \quad n=1,2, \ldots
\]
We refer to this sequence $\alpha: \mathbb{N}_0 \rightarrow
\mathbb{R}$ as the {\em deviation function} associated with the
scaling $P,K: \mathbb{N}_0 \rightarrow \mathbb{N}_0$. As the
terminology suggests, the deviation function measures by how much
the scaling deviates from the critical scaling $\frac{\log n}{n}$.

A scaling $P,K: \mathbb{N}_0 \rightarrow \mathbb{N}_0$ is said to
be {\em admissible} if
\begin{equation}
2 \leq  K_n \label{eq:AdmissibilityA}
\end{equation}
for {\em all} $n=1,2, \ldots $ {\em sufficiently} large. The main
result of this paper can now be stated as follows.

\begin{theorem}
{\sl Consider an admissible scaling $P,K: \mathbb{N}_0 \rightarrow
\mathbb{N}_0$ with deviation function $\alpha : \mathbb{N}_0
\rightarrow \mathbb{R}$ determined through
(\ref{eq:DeviationCondition}). We have
\[
\lim_{n \rightarrow \infty } P(n;\theta_n) =0 \quad \mbox{if~
$\lim_{ n\rightarrow \infty }\alpha_n = - \infty $.}
\label{eq:ZeroLaw+Connectivity}
\]
On the other hand, if there exists some $\sigma>0$ such that
\begin{equation}
\sigma n \leq P_n \label{eq:OneLaw+ConnectivityExtraCondition}
\end{equation}
for all $n=1,2, \ldots $ sufficiently large, then we have
\begin{equation}
\lim_{n \rightarrow \infty } P(n;\theta_n) = 1 \quad \mbox{if~
$\lim_{ n\rightarrow \infty }\alpha_n = \infty $.}
\label{eq:OneLaw+Connectivity}
\end{equation}
} \label{thm:ZeroOneLaw+Connectivity}
\end{theorem}

The condition (\ref{eq:OneLaw+ConnectivityExtraCondition}) is
sometimes expressed as $P_n = \Omega(n)$ and is slightly weaker
than the growth condition at (\ref{eq:OneBGMainCondition}) used by
Blackburn and Gerke \cite{BlackburnGerke}. Furthermore, Theorem
\ref{thm:ZeroOneLaw+Connectivity} implies the much weaker one-law
(\ref{eq:OneBGMainResult1})-(\ref{eq:OneBGMainResult2}). We also
note that the one-law in Theorem \ref{thm:ZeroOneLaw+Connectivity}
cannot hold if the condition (\ref{eq:AdmissibilityA}) fails. This
is a simple consequence of the following observation; see
\cite{YaganThesis} for details.

\begin{lemma}
{\sl For any mapping $P: \mathbb{N}_0 \rightarrow \mathbb{N}_0$
for which the limit $\lim_{n \rightarrow \infty} P_n$ exists
(possibly infinite), we have
\[
\lim_{n \rightarrow \infty } P(n;(1,P_n)) = \left \{
\begin{array}{ll}
0 & \mbox{if~ $\lim_{n \rightarrow \infty} P_n > 1$} \\
  &                 \\
1 & \mbox{if~ $\lim_{n \rightarrow \infty} P_n = 1$.}
\end{array}
\right . \label{eq:StrongZeroOneLaw+ConnectivityOnekey}
\]
} \label{lem:StrongZeroOneLaw+ConnectivityOnekey}
\end{lemma}

\section{A roadmap for the proof
of Theorem \ref{thm:ZeroOneLaw+Connectivity}} \label{sec:RoadMap}

Fix $n=2,3, \ldots $ and consider positive integers $K$ and $P$
such that $2 \leq K \leq P$. We define the events
\[
C_n (\theta) := \left [ \mathbb{K}(n; \theta) ~{\rm is~connected}
\right ]
\]
and
\[
I_n (\theta) := \left [ \mathbb{K}(n;\theta) ~{\rm
contains~no~isolated~nodes} \right ].
\]
If the random key graph $\mathbb{K}(n; \theta)$ is connected, then
it does not contain isolated nodes, whence $C_n(\theta)$ is a
subset of $I_n(\theta)$, and the conclusions
\begin{equation}
\bP{ C_n(\theta) } \leq \bP{ I_n(\theta) }
\label{eq:FromConnectivityToNodeIsolation1}
\end{equation}
and
\begin{equation}
\bP{ C_n(\theta)^c } = \bP{ C_n(\theta)^c \cap I_n (\theta) } +
\bP{ I_n(\theta)^c } \label{eq:FromConnectivityToNodeIsolation2}
\end{equation}
obtain.

In \cite{YaganMakowskiISIT2008}, we established the following
zero-one law for the absence of isolated nodes by the method of
first and second moments applied to the number of isolated nodes.

\begin{theorem}
{\sl For any admissible scaling $P,K: \mathbb{N}_0 \rightarrow
\mathbb{N}_0$, it holds that
\[
\lim_{n \rightarrow \infty } \bP{ I_n (\theta_n) } = \left \{
\begin{array}{ll}
0 & \mbox{if~ $\lim_{ n\rightarrow \infty }\alpha_n = - \infty $} \\
  &                 \\
1 & \mbox{if~ $\lim_{ n\rightarrow \infty }\alpha_n = + \infty $}
\end{array}
\right . \label{eq:VeryStrongZeroOneLaw+IsolatedNodes}
\]
where the deviation function $\alpha : \mathbb{N}_0 \rightarrow
\mathbb{R}$ is determined through (\ref{eq:DeviationCondition}). }
\label{thm:ZeroOneLaw+IsolatedNodes}
\end{theorem}

This result was also obtained independently by Blackburn and Gerke
\cite{BlackburnGerke} and Godehardt and Jaworski
\cite{GodehardtJaworski}. In this last paper the authors show the
stronger result that the number of isolated nodes is
asymptotically Poisson distributed with parameter $e^{-c}$ under
scalings of the form (\ref{eq:DeviationCondition}) with deviation
function satisfying $\lim_{n \rightarrow \infty } \alpha_n = c$
for some finite scalar $c$.

Taken together with Theorem \ref{thm:ZeroOneLaw+IsolatedNodes},
the relations (\ref{eq:FromConnectivityToNodeIsolation1}) and
(\ref{eq:FromConnectivityToNodeIsolation2}) pave the way to
proving Theorem \ref{thm:ZeroOneLaw+Connectivity}. Indeed, pick an
admissible scaling $P,K: \mathbb{N}_0 \rightarrow \mathbb{N}_0$
with deviation function $\alpha : \mathbb{N}_0 \rightarrow
\mathbb{R}$. If $\lim_{n \rightarrow \infty} \alpha_n = -\infty$,
then $\lim_{n \rightarrow \infty} \bP{ I_n(\theta_n) } = 0$ by the
zero-law for the absence of isolated nodes, whence $\lim_{n
\rightarrow \infty} \bP{ C_n(\theta_n) } = 0$ with the help of
(\ref{eq:FromConnectivityToNodeIsolation1}). If $\lim_{n
\rightarrow \infty} \alpha_n = \infty$, then $\lim_{n \rightarrow
\infty} \bP{ I_n(\theta_n) } = 1$ by the one-law for the absence
of isolated nodes, and the desired conclusion $\lim_{n \rightarrow
\infty } \bP{ C_n(\theta_n) } = 1$ (or equivalently, $\lim_{n
\rightarrow \infty } \bP{ C_n(\theta_n)^c } = 0$) will follow via
(\ref{eq:FromConnectivityToNodeIsolation2}) if we show that
\begin{equation}
\lim_{n \rightarrow \infty } \bP{ C_n(\theta_n)^c \cap I_n
(\theta_n) } = 0. \label{eq:RoadMap1}
\end{equation}
We shall do this by finding a sufficiently tight upper bound on
the probability in (\ref{eq:RoadMap1}) and then showing that it
goes to zero as well. While the additional condition
(\ref{eq:OneLaw+ConnectivityExtraCondition}) plays a crucial role
in carrying out this argument, a number of additional assumptions
will be imposed on the admissible scaling under consideration.
This is done mostly for technical reasons in that it leads to
simpler proofs. Eventually these additional conditions will be
removed to ensure the desired final result, namely
(\ref{eq:OneLaw+Connectivity}) under
(\ref{eq:OneLaw+ConnectivityExtraCondition}), e.g., see Section
\ref{sec:ReductionStep} for details.

With this in mind, the admissible scaling $P,K: \mathbb{N}_0
\rightarrow \mathbb{N}_0$ is said to be {\em strongly admissible}
if its deviation function $\alpha: \mathbb{N}_0 \rightarrow
\mathbb{R}$ satisfies the additional growth condition
\begin{equation}
\alpha_n = o(n). \label{eq:AdmissibilityB}
\end{equation}

Strong admissibility has the following useful implications: Under
(\ref{eq:AdmissibilityB}) it is always the case from
(\ref{eq:DeviationCondition}) that
\begin{equation}
\lim_{n \rightarrow \infty} \frac{K^2_n}{P_n} = 0.
\label{eq:RatioConditionStrong+Consequence0}
\end{equation}
Since $1 \leq K_n \leq K_n^2$ for all $n=1,2, \ldots $, this last
convergence implies
\begin{equation}
\lim_{n \rightarrow \infty } \frac{K_n}{P_n} = 0 \quad \mbox{and}
\quad \lim_{n \rightarrow \infty } P_n = \infty .
\label{eq:RatioConditionStrong+Consequence1}
\end{equation}
As a result, we have
\begin{equation}
2 K_n \leq P_n \label{eq:RatioConditionStrong+Consequence3}
\end{equation}
for all $n=1,2, \ldots $ sufficiently large, and the random key
graph does not degenerate into a complete graph under a strongly
admissible scaling. Finally, in Lemma
\ref{lem:AsymptoticEquivalence} we show that
(\ref{eq:RatioConditionStrong+Consequence0}) suffices to imply
\begin{equation}
1 - q(\theta_n) \sim \frac{K^2_n}{P_n}.
\label{eq:EquivalentCondition}
\end{equation}
This is discussed in Section \ref{sec:Usefulequivalence}, and
provides the appropriate version of (\ref{eq:EquivalenceNeeded}).

\section{A reduction step}
\label{sec:ReductionStep}

The relevance of the notion of strong admissibility flows from the
following fact.

\begin{lemma}
{\sl Consider an admissible scaling $K,P: \mathbb{N}_0 \rightarrow
\mathbb{N}_0$ whose deviation sequence $\alpha : \mathbb{N}_0
\rightarrow \mathbb{R}$ satisfies
\[
\lim_{n \rightarrow \infty} \alpha_n = \infty.
\label{eq:ReductionStep2}
\]
Assume there exists some $\sigma > 0$ such that
(\ref{eq:OneLaw+ConnectivityExtraCondition}) holds for all $n=1,2,
\ldots $ sufficiently large. Then, there always exists an
admissible scaling $\tilde K,\tilde P: \mathbb{N}_0 \rightarrow
\mathbb{N}_0$ with
\begin{equation}
\tilde K_n \leq K_n \quad {\rm and} \quad \tilde P_n = P_n, \quad
n=1,2, \ldots \label{eq:ComparingAlpha+A}
\end{equation}
whose deviation function $\tilde \alpha : \mathbb{N}_0 \rightarrow
\mathbb{R}$ satisfies both conditions
\begin{equation}
\lim_{n \rightarrow \infty} \tilde \alpha_n = \infty \quad
\mbox{and} \quad \tilde \alpha_n = o (n).
\label{eq:ComparingAlpha+B}
\end{equation}
} \label{lem:ReductionStep}
\end{lemma}

\myproof For each $n=1,2, \ldots $, set
\[
K^\star_n := \sqrt{ P_n \cdot \frac{ \log n + \alpha^\star_n }{n}
} \quad \mbox{where} \ \alpha^\star_n := \min \left ( \alpha_n ,
\log n \right ).
\]
The properties
\begin{equation}
\lim_{n \rightarrow \infty} \alpha^\star_n = \infty
\label{eq:GrowthConditionForAlphaStar1}
\end{equation}
and
\begin{equation}
\alpha^\star_n = o(n) \label{eq:GrowthConditionForAlphaStar2}
\end{equation}
are immediate by construction.

Now define the scaling $\tilde K, \tilde P: \mathbb{N}_0
\rightarrow \mathbb{N}_0$ by
\[
\tilde K_n := \left \lceil K^\star_n \right \rceil , \ \tilde P_n
= P_n , \quad n=1,2, \ldots \label{eq:ComparingAlpha+New}
\]
We get $K^\star_n \leq K_n$ for all $n=1,2, \ldots $ since
$\alpha^\star_n \leq \alpha_n$, whence $\tilde K_n \leq K_n$ by
virtue of the fact that $K_n$ is always an integer. This
establishes (\ref{eq:ComparingAlpha+A}).

Next, observe that $\tilde K_n = 1 $ if and only $K^\star_n \leq 1
$, a condition which occurs only when
\begin{equation}
P_n \left ( \log n + \alpha^\star_n \right ) \leq n.
\label{eq:ConditionOnAFiniteSubset}
\end{equation}
This last inequality can only hold for a finite number of values
of $n$. Otherwise, there would exist a {\em countably infinite}
subset $N$ of $\mathbb{N}_0$ such that both
(\ref{eq:OneLaw+ConnectivityExtraCondition}) and
(\ref{eq:ConditionOnAFiniteSubset}) simultaneously hold on $N$. In
that case, we conclude that
\[
\sigma \left ( \log n + \alpha^\star_n \right ) \leq 1, \quad n
\in N
\]
and this is a clear impossibility in view of
(\ref{eq:GrowthConditionForAlphaStar1}). Together with
(\ref{eq:ComparingAlpha+A}) this establishes the admissibility of
the scaling $\tilde K, \tilde P: \mathbb{N}_0 \rightarrow
\mathbb{N}_0$.

Fix $n=1,2, \ldots $. The definitions imply $K^\star_n \leq \tilde
K_n < 1 + K^\star_n$, and upon squaring we get the inequalities
\[
P_n \cdot \frac{ \log n + \alpha^\star_n }{n} \leq \tilde K^2_n
\]
and
\[
\tilde K^2_n < 1 + 2 \sqrt{ P_n \cdot \frac{ \log n +
\alpha^\star_n }{n} } + P_n \cdot \frac{ \log n + \alpha^\star_n
}{n}.
\]

The deviation sequence $\tilde \alpha: \mathbb{N}_0 \rightarrow
\mathbb{R}$ of the newly defined scaling
(\ref{eq:ComparingAlpha+A}) is determined through
\[
\frac{ \tilde K^2_n}{\tilde P_n} = \frac{\log n + \tilde
\alpha_n}{n}, \quad n=1,2, \ldots .
\]
Using the two inequalities above we then conclude that
\begin{equation}
\alpha^\star_n \leq \tilde \alpha_n
\label{eq:EquivalentCondition1}
\end{equation}
and
\begin{equation}
\frac{ \tilde \alpha_n }{n} < \frac{1}{P_n} + 2 \sqrt{
\frac{1}{P_n} \cdot \frac{ \log n + \alpha^\star_n }{n} } + \frac{
\alpha^\star_n }{n}. \label{eq:EquivalentCondition1bis}
\end{equation}
It is now plain from (\ref{eq:GrowthConditionForAlphaStar1}) and
(\ref{eq:EquivalentCondition1}) that the first half of
(\ref{eq:ComparingAlpha+B}) holds. Next, by combining
(\ref{eq:EquivalentCondition1}) and
(\ref{eq:EquivalentCondition1bis}) we get
\begin{equation}
\frac{\alpha^\star_n}{n} \leq \frac{ \tilde \alpha_n }{n} <
\frac{1}{P_n} + 2 \sqrt{ \frac{1}{P_n} \cdot \frac{ \log n +
\alpha^\star_n }{n} } + \frac{ \alpha^\star_n }{n}.
\label{eq:EquivalentCondition2}
\end{equation}
Letting $n$ go to infinity in (\ref{eq:EquivalentCondition2}) and
using (\ref{eq:GrowthConditionForAlphaStar2}) we conclude to the
second half of (\ref{eq:ComparingAlpha+B}) since
$\lim_{n\rightarrow \infty} P_n = \infty$ by virtue of
(\ref{eq:OneLaw+ConnectivityExtraCondition}). \myendpf


The scaling $\tilde K,\tilde P: \mathbb{N}_0 \rightarrow
\mathbb{N}_0$ defined at (\ref{eq:ComparingAlpha+A}) is strongly
admissible and still satisfies the condition
(\ref{eq:OneLaw+ConnectivityExtraCondition}), and an easy coupling
argument based on (\ref{eq:ComparingAlpha+A}) shows that
\[
P(n; \tilde \theta_n ) \leq P(n; \theta_n ), \quad n=2,3, \ldots
\]
Therefore, we need only show (\ref{eq:OneLaw+Connectivity}) under
(\ref{eq:OneLaw+ConnectivityExtraCondition}) for strongly
admissible scalings. As a result, in view of the discussion
leading to (\ref{eq:RoadMap1}) it suffices to establish the
following result, to which the remainder of the paper is devoted.

\begin{proposition}
{\sl Consider any strongly admissible scaling $P,K: \mathbb{N}_0
\rightarrow \mathbb{N}_0$ whose deviation function $\alpha:
\mathbb{N}_0 \rightarrow \mathbb{R}$ satisfies $\lim_{n
\rightarrow \infty} \alpha_n = \infty$. Under the condition
(\ref{eq:OneLaw+ConnectivityExtraCondition}), we have
\begin{equation}
\lim_{n \rightarrow \infty} \bP{ C_n (\theta_n)^c \cap I_n
(\theta_n) } = 0 . \label{eq:OneLawAfterReduction}
\end{equation}
} \label{prop:OneLawAfterReduction}
\end{proposition}

Proposition \ref{prop:OneLawAfterReduction} shows that in random
key graphs, graph connectivity is asymptotically equivalent to the
absence of isolated nodes under any strongly admissible scaling
whose deviation function $\alpha: \mathbb{N}_0 \rightarrow
\mathbb{R}$ satisfies $\lim_{n \rightarrow \infty} \alpha_n =
\infty$ under the condition
(\ref{eq:OneLaw+ConnectivityExtraCondition}).

\section{The equivalence (\ref{eq:EquivalentCondition})}
\label{sec:Usefulequivalence}

To establish the key equivalence (\ref{eq:EquivalentCondition}) we
start with simple bounds which prove useful in a number of places.
Full details are available in
\cite{YaganMakowskiISRTechReport2009, YaganThesis}.

\begin{lemma}
{\sl For positive integers $K$, $L$ and $P$ such that $K + L \leq
P$, we have
\[
\left ( 1 - \frac{L}{P-K} \right )^K \leq \frac{{P- L \choose
K}}{{P \choose K}} \leq \left ( 1 - \frac{L}{P} \right )^K,
\label{eq:BoundsOnRatioCombinatorial}
\]
whence
\begin{equation}
\frac{{P- L \choose K}}{{P \choose K}} \leq e^{-K \cdot
\frac{L}{P}} . \label{eq:BoundsOnRatioCombinatorial+2}
\end{equation}
} \label{lem:BoundsOnRatioCombinatorial}
\end{lemma}

Applying Lemma \ref{lem:BoundsOnRatioCombinatorial} (with $L=K$)
to the expression (\ref{eq:q_theta}) yields the following bounds.

\begin{lemma}
{\sl With positive integers $K$ and $P$ such that $2K \leq P$, we
have
\[
1 - e^{-\frac{K^2}{P}} \leq 1 - q(\theta) \leq \frac{K^2}{P-K} .
\label{eq:BoundsOnOneMinusQ(Theta)1}
\]
} \label{lem:BoundsOnOneMinusQ(Theta)}
\end{lemma}
A little bit more than (\ref{eq:EquivalentCondition}) can then be
said.

\begin{lemma}
{\sl For any scaling $P,K: \mathbb{N}_0 \rightarrow \mathbb{N}_0$,
it holds that
\begin{equation}
\lim_{n \rightarrow \infty} q(\theta_n) = 1 \label{eq:Condition1}
\end{equation}
if and only if
\begin{equation}
\lim_{n \rightarrow \infty} \frac{K^2_n}{P_n} = 0,
\label{eq:Condition2}
\end{equation}
and under either condition we have the asymptotic equivalence
\begin{equation}
1 - q(\theta_n) \sim \frac{K^2_n}{P_n}.
\label{eq:AsymptoticsEquivalence}
\end{equation}
} \label{lem:AsymptoticEquivalence}
\end{lemma}

On several occasions, we will rely on
(\ref{eq:AsymptoticsEquivalence}) through the following equivalent
formulation: For every $\delta$ in $(0,1)$ there exists a finite
integer $n^\star (\delta)$ such that
\begin{equation}
(1-\delta) \frac{K_n^2}{P_n} \leq 1 - q (\theta_n) \leq (1+\delta)
\frac{K_n^2}{P_n} \label{eq:AsymptoticEquivalenceWithDelta}
\end{equation}
whenever $n \geq n^\star(\delta)$.

\myproof As noted already at the end of Section \ref{sec:RoadMap},
condition (\ref{eq:Condition2}) (which holds for any strongly
admissible scaling) implies
(\ref{eq:RatioConditionStrong+Consequence3}) for all $n=1,2,
\ldots$ sufficiently large. On that range Lemma
\ref{lem:BoundsOnOneMinusQ(Theta)} yields
\begin{equation}
1 - e^{-\frac{K^2_n}{P_n}} \leq 1 - q (\theta_n) \leq
\frac{K^2_n}{P_n-K_n} . \label{eq:Equivalence1}
\end{equation}

Multiply (\ref{eq:Equivalence1}) by $\frac{P_n}{K^2_n}$ and let
$n$ go to infinity in the resulting set of inequalities. Under
(\ref{eq:Condition2}), we get
\[
\lim_{n\rightarrow\infty} \frac{P_n}{K^2_n} \cdot \left ( 1 -
e^{-\frac{K^2_n}{P_n}} \right ) = 1 \label{eq:EquivalenceLimit1}
\]
from the elementary fact $\lim _{t \downarrow 0} \frac{ 1 - e^{-t}
}{t} = 1$, while
\[
\lim_{n\rightarrow\infty} \frac{P_n}{K^2_n} \cdot \frac{K^2_n}{P_n
- K_n} = \lim_{n\rightarrow\infty} \frac{P_n}{P_n - K_n } =1
\label{eq:EquivalenceLimit2}
\]
by virtue of (\ref{eq:RatioConditionStrong+Consequence1}) (which
is implied by (\ref{eq:Condition2})). The asymptotic equivalence
(\ref{eq:AsymptoticsEquivalence}) follows, and the validity of
(\ref{eq:Condition1}) is immediate.

Conversely, under the condition $\lim_{n \rightarrow \infty}
q(\theta_n) =1 $, we have $0 < q(\theta_n) < 1 $ for all $n$
sufficiently large (by the comment following
(\ref{eq:Probab_key_ring_does_not_intersect_S})), and the
constraint (\ref{eq:RatioConditionStrong+Consequence3})
necessarily holds for all $n=1,2, \ldots$ sufficiently large. On
that range, (\ref{eq:Equivalence1}) being valid, we conclude to
$\lim_{n \rightarrow \infty}e^{-\frac{K^2_n}{P_n}} =1$ under
(\ref{eq:Condition1}). The convergence (\ref{eq:Condition2}) now
follows and the asymptotic equivalence
(\ref{eq:AsymptoticsEquivalence}) is given by the first part of
the proof. \myendpf

\section{A basic union bound}
\label{sec:BasicUnionBound}

Proposition \ref{prop:OneLawAfterReduction} will be established
with the help of a union bound for the probability appearing at
(\ref{eq:OneLawAfterReduction}) -- The approach is similar to the
one used for proving the one-law for connectivity in
Erd\H{o}s-R\'enyi graphs \cite[p. 164]{Bollobas}
\cite[Section 3.4, p. 40]{DraiefMassoulieLMS} 
\cite[p. 304]{SpencerSaintFlour1991}:

Fix $n=2,3, \ldots $ and consider positive integers $K$ and $P$
such that $2K \leq P$. For any non-empty subset $S$ of nodes,
i.e., $S \subseteq \{1, \ldots , n \}$, we define the graph
$\mathbb{K} (n;\theta) (S)$ (with vertex set $S$) as the subgraph
of $\mathbb{K} (n;\theta)$ restricted to the nodes in $S$. We also
say that $S$ is {\em isolated} in $\mathbb{K} (n;\theta)$ if there
are no edges (in $\mathbb{K} (n;\theta)$) between the nodes in $S$
and the nodes in the complement $S^c = \{ 1, \ldots , n \} - S$.
This is characterized by
\[
K_i (\theta) \cap K_j (\theta) = \emptyset , \quad i \in S , \ j
\in S^c  .
\]

With each non-empty subset $S$ of nodes, we associate several
events of interest: Let $C_n (\theta ; S)$ denote the event that
the subgraph $\mathbb{K} (n;\theta) (S)$ is itself connected. The
event $C_n (\theta ; S)$ is completely determined by the rvs $\{
K_i(\theta), \ i \in S \}$. We also introduce the event $B_n
(\theta ; S)$ to capture the fact that $S$ is isolated in
$\mathbb{K} (n;\theta)$, i.e.,
\[
B_n (\theta ; S) := \left [ K_i (\theta) \cap K_j (\theta) =
\emptyset , \quad i \in S , \ j \in S^c \right ] .
\]
Finally, we set
\begin{equation}
A_n (\theta ; S) := C_n (\theta ; S) \cap B_n (\theta ; S) .
\label{eq:A_{n,r}(theta)}
\end{equation}

The starting point of the discussion is the following basic
observation: If $\mathbb{K} (n;\theta)$ is {\em not} connected and
yet has {\em no} isolated nodes, then there must exist a subset
$S$ of nodes with $|S| \geq 2$ such that $\mathbb{K}(n;\theta)
(S)$ is connected while $S$ is isolated in $\mathbb{K}
(n;\theta)$. This is captured by the inclusion
\[
C_n(\theta)^c \cap I_n(\theta) \subseteq \cup_{S \in \mathcal{N}:
~ |S| \geq 2} ~ A_n (\theta ; S) \label{eq:BasicIdea}
\]
with $\mathcal{N}$ denoting the collection of all non-empty
subsets of $\{ 1, \ldots , n \}$. A moment of reflection should
convince the reader that this union need only be taken over all
subsets $S$ of $\{1, \ldots , n \}$ with $2 \leq |S| \leq \lfloor
\frac{n}{2} \rfloor $. Then, a standard union bound argument
immediately gives
\begin{eqnarray}
\bP{ C_n(\theta)^c \cap I_n(\theta) } &\leq & \sum_{ S \in
\mathcal{N}: 2 \leq |S| \leq \lfloor \frac{n}{2} \rfloor } \bP{
A_n (\theta ; S) }
\nonumber \\
&=& \sum_{r=2}^{ \lfloor \frac{n}{2} \rfloor } \left ( \sum_{S \in
\mathcal{N}_r } \bP{ A_n (\theta ; S) } \right )
\label{eq:BasicIdea+UnionBound}
\end{eqnarray}
where $\mathcal{N}_{r} $ denotes the collection of all subsets of
$\{ 1, \ldots , n \}$ with exactly $r$ elements.

For each $r=1, \ldots , n$, we simplify the notation by writing
$A_{n,r} (\theta) := A_n (\theta ; \{ 1, \ldots , r \} )$,
$B_{n,r} (\theta) := B_n (\theta ; \{ 1, \ldots , r \} )$ and
$C_{r} (\theta) := C_n (\theta ; \{ 1, \ldots , r \} )$. For $r=n$
this notation is consistent with $C_n(\theta)$ as defined in
Section \ref{sec:RoadMap}. Under the enforced assumptions,
exchangeability gives
\[
\bP{ A_n (\theta ; S) } = \bP{ A_{n,r} (\theta ) }, \quad S \in
\mathcal{N}_r
\]
and the expression
\[
\sum_{S \in \mathcal{N}_r } \bP{ A_n (\theta ; S) } = {n \choose
r} ~ \bP{ A_{n,r}(\theta ) } \label{eq:ForEach=r}
\]
follows since $|\mathcal{N}_{r} | = {n \choose r}$. Substituting
into (\ref{eq:BasicIdea+UnionBound}) we obtain the key bound
\begin{equation}
\bP{ C_n(\theta)^c \cap I_n(\theta) } \leq \sum_{r=2}^{ \lfloor
\frac{n}{2} \rfloor } {n \choose r} ~ \bP{ A_{n,r}(\theta ) } .
\label{eq:BasicIdea+UnionBound2}
\end{equation}

Consider a strongly  admissible scaling $P,K: \mathbb{N}_0
\rightarrow \mathbb{N}_0$ as in the statement of Proposition
\ref{prop:OneLawAfterReduction}. In the right hand side of
(\ref{eq:BasicIdea+UnionBound2}) we substitute $\theta$ by
$\theta_n$ by means of this strongly admissible scaling. The proof
of Proposition \ref{prop:OneLawAfterReduction} will be completed
once we show that
\begin{equation}
\lim_{n \rightarrow \infty} \sum_{r=2}^{ \lfloor \frac{n}{2}
\rfloor } {n \choose r} ~ \bP{ A_{n,r}(\theta_n) } = 0
\label{eq:OneLawToShow}
\end{equation}
under the appropriate conditions. This approach was used to
establish the one-law in Erd\H{o}s-R\'enyi graphs \cite{Bollobas,
DraiefMassoulieLMS, SpencerSaintFlour1991} where simple bounds can
be derived for the probability terms in (\ref{eq:OneLawToShow}).
Our situation is technically more involved and requires more
delicate bounding arguments as will become apparent in the
forthcoming sections.

\section{Bounding the probabilities $\bP{A_{n,r}(\theta)}$ \\
         ($r=1, \ldots , n$)}
\label{sec:BoundingProbabilities}

Again consider positive integers $K$ and $P$ such that $2K \leq
P$. Fix $n=2,3, \ldots $ and pick $r=1, \ldots , n-1$. Since exact
expressions are not available for the probability $\bP{
A_{n,r}(\theta) }$, we seek instead to provide a bound on this
quantity. For reasons that will become apparent shortly, it will
be beneficial to focus on the following more general task: Let
${\cal F}_r$ denotes the $\sigma$-field on $\Omega$ generated by
the rvs $K_1(\theta), \ldots , K_r(\theta)$. We
are interested in deriving an upper bound on the probability $\bP{
A_{n,r}(\theta)  \cap E }$ where $E$ is {\em any} ${\cal
F}_r$-measurable event, the original situation corresponding to $E
= \Omega$.

In the course of doing so, we shall make use of the rv $U_{r}
(\theta)$ given by
\[
U_{r} (\theta) := \left | \cup_{i=1}^r K_i(\theta) \right | .
\label{eq:U}
\]
The rv $U_r(\theta)$ counts the number of {\em distinct} keys
issued to the nodes $1, \ldots , r$, so that the bounds
\begin{equation}
K \leq U_{r} (\theta) \leq \min \left ( rK,P \right )
\label{eq:BoundOnU}
\end{equation}
always hold.

Thus, pick any ${\cal F}_r$-measurable event $E$, and note that
$C_r(\theta)$ is also an ${\cal F}_r$-measurable event since
completely determined by the rvs $K_1(\theta), \ldots ,
K_r(\theta) $. It is now plain (\ref{eq:A_{n,r}(theta)}) that
\begin{eqnarray}
\bP{ A_{n,r}(\theta)  \cap E } &=& \bP{ B_{n,r}(\theta) \cap
C_r(\theta) \cap E }
\nonumber \\
&=& \bE{ \1{ C_r(\theta) \cap E } \bP{  B_{n,r}(\theta) | {\cal
F}_r } } \nonumber
\end{eqnarray}
upon preconditioning on the rvs $K_1(\theta), \ldots , K_r(\theta)
$. Next, with the help of the equivalence
\[
B_{n,r}(\theta ) = \left [ \left ( \cup_{i=1}^r K_i(\theta)
\right) \cap K_j(\theta) = \emptyset, \ j=r+1, \ldots n \right ],
\]
we can use (\ref{eq:Probab_key_ring_does_not_intersect_S}) (with
$S = \cup_{i=1}^r K_i(\theta) $) to get
\begin{eqnarray}
\lefteqn{ \bP{ B_{n,r} (\theta) | {\cal F}_r } } & &
\label{eq:ExpressionConditionalProbability} \nonumber \\
&=& \left ( \frac{{P- U_{r} (\theta) \choose K}}{{P \choose K}}
\right )^{n-r} \1{ U_r(\theta) \leq P - K } \quad a.s. \nonumber
\end{eqnarray}
under the enforced independence assumptions. The conclusion
\begin{eqnarray}
\bP{ A_{n,r}(\theta ) \cap E } = \bE{ \1{ C^\star_{r}(\theta )
\cap E } \cdot \left ( \frac{{P- U_{r} (\theta) \choose K}}{{P
\choose K}} \right )^{n-r} } \label{eq:ComputePA_{n,r}} \nonumber
\end{eqnarray}
then follows with
\[
C^\star_{r}(\theta) := C_{r}(\theta) \cap [ U_{r} (\theta) \leq P
- K ] .
\]
Applying (\ref{eq:BoundsOnRatioCombinatorial+2}) (with $L =
U_r(\theta)$) in Lemma \ref{lem:BoundsOnRatioCombinatorial}, we
finally obtain the inequality
\begin{eqnarray}
\lefteqn{ \bP{ A_{n,r}(\theta) \cap E }  } & &
\nonumber \\
&\leq& \bE{ \1{ C^\star_{r}(\theta) \cap E } \cdot e^{-(n-r)
\frac{K}{P} \cdot U_{r} (\theta) } } .
\label{eq:BoundComputePA_{n,r}}
\end{eqnarray}

This discussion already brings out a number of items that are
likely to require some attention: We will need good bounds for the
probabilities $\bP{ C_{r} (\theta) }$ and $\bP{ C^\star _{r}
(\theta) }$. Also, some of the distributional properties of the rv
$U_{r}(\theta)$ are expected to play a role. The constraints
(\ref{eq:BoundOnU}) automatically imply $U_{r} (\theta) \leq P - K
$ whenever $rK \leq P-K$, i.e., $(r+1) K \leq P$, whence
\begin{equation}
C^\star_{r}(\theta) = C_{r}(\theta), \quad r=1, \ldots ,
r_n(\theta)
\end{equation}
where we have set
\[
r_n (\theta) := \min \left ( r(\theta), \left \lfloor \frac{n}{2}
\right \rfloor \right ) \quad {\rm with} \quad r(\theta) := \left
\lfloor \frac{P}{K} \right \rfloor - 1.
\]
This suggests that different arguments will probably be needed for
the ranges $1\leq r \leq r_n(\theta)$ and $r_n(\theta) < r \leq
\lfloor \frac{n}{2} \rfloor$.

The next result is crucial to showing that for each $r=2. \ldots ,
n$, the probability of the event $C_{r} (\theta)$ can be provided
an upper bound in terms of known quantities. Let
$\mathbb{K}_r(n;\theta)$ stand for the subgraph
$\mathbb{K}(n;\theta)(S)$ when $S = \{1 , \ldots , r \}$, and let
${\cal T}_{r}$ denote the collection of all spanning trees on the
vertex set $\{1, \ldots , r \}$.

\begin{lemma}
{\sl For each $r=2, \ldots , n$, we have
\begin{equation}
\bP{ T \subset \mathbb{K}_r(n;\theta) } = \left ( 1 - q(\theta)
\right )^{r-1}, \quad T \in {\cal T}_{r}
\label{eq:ProbabilityOfTree}
\end{equation}
where the notation $T \subset \mathbb{K}_r(n;\theta)$ indicates
that the tree $T$ is a subgraph spanning $\mathbb{K}_r(n;\theta)$.
} \label{lem:ProbabilityOfTree}
\end{lemma}

This last expression is analogous to the one found in
Erd\H{o}s-R\'enyi graphs \cite{Bollobas, DraiefMassoulieLMS} 
with $1-q(\theta)$
playing the role of probability of link assignment, and this in
spite of the correlations between some link assignments.

\myproof We shall prove the result by induction on $r=2, \ldots ,
n$. For $r=2$ the conclusion (\ref{eq:ProbabilityOfTree}) is
nothing more than (\ref{eq:q_theta}) since ${\cal T}_2$ contains
exactly one tree, and this establishes the basis step.

Next, we consider the following induction step: Pick $r=2, \ldots
, n-1$ and assume that for each $s=2, \ldots , r$, it is already
known that
\begin{equation}
\bP{ T \subset \mathbb{K}_s(n;\theta) } = \left ( 1 - q(\theta)
\right )^{s-1}, \quad T \in {\cal T}_{s} .
\label{eq:ProbabilityOfTreeInduction}
\end{equation}
We now show that (\ref{eq:ProbabilityOfTreeInduction}) also holds
for each $s=2, \ldots , ,r+1$. To that end, pick a tree $T$ in
${\cal T}_{r+1}$ and identify its root.\footnote{As we are
considering undirected graphs, all nodes can act as a root for the
(undirected) tree $T$, in which case any one will do for the
forthcoming discussion.} Let $i$ denote a node that is farthest
from the root of $T$ -- There might be several such nodes. Also
denote by $p$ its unique parent, and let $D(p)$ denote the set of
children of $p$. Obviously $D(p)$ is not empty as it contains node
$i$; set $|D(p)| = d$. Next we construct a new tree $T^\star$ from
$T$ by removing from $T$ all the edges from node $p$ to the nodes
in $D(p)$. By exchangeability, there is no loss of generality in
assuming (as we do from now on) that the tree is rooted at node
$1$, that the unique parent $p$ of the farthest node selected has
label $r-d+1$, and that its children have been labelled $r-d+2,
\ldots ,r+1$. With this convention, the tree $T^\star$ is defined
on the set of nodes $\{1, \ldots , r-d+1 \}$.

It is plain that $T \subseteq \mathbb{K}_{r+1} (n,; \theta ) $
occurs if and only if the two sets of conditions
\[
K_{r-d+1} (\theta ) \cap K_\ell (\theta ) \neq \emptyset, \quad
\ell = r-d+2, \ldots ,r+1
\]
and
\[
T^\star \subseteq \mathbb{K}_{r-d+1} (n; \theta )
\]
both hold. Under the enforced independence assumptions we get
\[
\bP{
\begin{array}{c}
K_{r-d+1} (\theta ) \cap K_\ell (\theta ) \neq \emptyset, \\
\ell = r-d+2, \ldots ,r+1 \\
\end{array}
\Big | {\cal F}_{r-d+1} } = (1-q(\theta))^{d}.
\]
Thus, upon conditioning with respect to the rvs $K_1 (\theta),
\ldots , K_{r-d+1}(\theta) $ we readily find
\begin{eqnarray}
\lefteqn{ \bP{ T \subseteq \mathbb{K}_{r+1} (n,; \theta ) } } & &
\nonumber \\
&=& (1-q(\theta))^{d} \bP{ T^\star \subseteq \mathbb{K}_{r-d+1}
(n; \theta ) }
\nonumber \\
&=& (1-q(\theta) )^{d} (1-q(\theta) )^{r- d}
\nonumber \\
&=& (1-q(\theta) )^{r} \nonumber \nonumber
\end{eqnarray}
as we use the induction hypothesis
(\ref{eq:ProbabilityOfTreeInduction}) to evaluate the probability
of the event $[T^\star \subseteq \mathbb{K}_{r-d+1} (n; \theta)]$.
This establishes the induction step. \myendpf

The bound below now follows as in Erd\H{o}s-R\'enyi graphs
\cite{Bollobas,DraiefMassoulieLMS}.

\begin{lemma}
{\sl For each $r=2, \ldots , n$, we have
\begin{equation}
\bP{ C_{r} (\theta) } \leq r^{r-2} \left ( 1 - q(\theta)
\right)^{r-1} . \label{eq:ProbabilityOfC}
\end{equation}
} \label{lem:ProbabilityOfC}
\end{lemma}

\myproof Fix $r=2, \ldots , n$. If $\mathbb{K}_r (n;\theta)$ is a
connected graph, then it must contain a spanning tree on the
vertex set $\{1, \ldots . r \}$, and a union bound argument yields
\[
\bP{ C_{r} (\theta) } \leq \sum_{T \in {\cal T}_r } \bP{ T \subset
\mathbb{K}(n;\theta)(S) } . \label{eq:ProbabilityOfAwithTrees}
\]
By Cayley's formula \cite{CayleysFormula, Martin} there are
$r^{r-2}$ trees on $r$ vertices, i.e., $| {\cal T}_{r}| =
r^{r-2}$, and (\ref{eq:ProbabilityOfC}) follows upon making use of
(\ref{eq:ProbabilityOfTree}). \myendpf

The bound (\ref{eq:BoundComputePA_{n,r}}) (with $E=\Omega$) and
the inequality $U_r(\theta) \geq K$ together imply
\begin{eqnarray}
\bP{ A_{n,r}(\theta ) } &\leq& \bP{ C_{r}(\theta ) } \cdot
e^{-(n-r) \frac{K^2}{P} }
\nonumber \\
&\leq & r^{r-2} \left ( 1 - q(\theta) \right )^{r-1} \cdot
e^{-(n-r) \frac{K^2}{P} } \label{eq:CrudeUpperBound}
\end{eqnarray}
as we make use of Lemma \ref{lem:ProbabilityOfC} in the last step.
Unfortunately, this bound turns out to be too loose for our
purpose. As this can be traced to the crude lower bound used for
$U_r(\theta)$, we expect that improvements are possible if we take
into account the distributional properties of the rv
$U_r(\theta)$. This step is taken in the next section.

\section{The tail of the rv $U_{r}(\theta)$ and improved bounds}
\label{sec:Tail+ImprovedBounds}

Consider positive integers $K$ and $P$ such that $K \leq P$. Rough
estimates will suffice to get the needed information regarding the
distribution of the rv $U_r(\theta)$. This is the content of the
next result.

\begin{lemma}
{\sl For all $r=1,2, \ldots $, the bounds
\begin{equation}
\bP{ U_{r}(\theta) \leq x } \leq {P \choose x} \left ( \frac{x}{P}
\right )^{rK} \label{eq:BoundOnProbability+2}
\end{equation}
holds whenever $x=K, \ldots , \min ( rK, P )$. }
\label{lem:BoundOnProbability+2}
\end{lemma}

\myproof For a given $x$ in the prescribed range, we note that
$U_{r}(\theta) \leq x $ implies that $\cup_{i=1}^r K_i(\theta) $
is contained in some set $S$ of size $x$, whence
\[
[ U_{r}(\theta) \leq x ] \subseteq \bigcup_{S \in \mathcal{P}_x} [
\cup_{i=1}^r K_i(\theta) \subseteq  S ].
\]
A standard union bound argument gives
\begin{eqnarray}
\bP{ U_{r}(\theta) \leq x } &\leq & \sum_{S \in \mathcal{P}_x}
\bP{ \cup_{i=1}^r K_i(\theta) \subseteq S }
\nonumber \\
&=& \sum_{S \in \mathcal{P}_x} \bP{ K_i(\theta) \subseteq S, \
i=1, \ldots , r }
\nonumber \\
&=& \sum_{S \in \mathcal{P}_x} \prod_{i=1}^r \bP{ K_i(\theta)
\subseteq S }
\nonumber \\
&=& \sum_{S \in \mathcal{P}_x} \left ( \bP{ K_1(\theta) \subseteq
S } \right )^r \label{eq:BoundOnProbability+1A}
\end{eqnarray}
under the enforced assumptions on the rvs $K_1(\theta), \ldots ,
K_n(\theta) $.

Since every subset of size $x$ contains ${x \choose K}$ further
subsets of size $K$, we get
\[
\bP{ K_1(\theta) \subseteq S } = \frac{ {x \choose K } }{ {P
\choose K} }, \quad S \in \mathcal{P}_x .
\]
Substituting this fact into (\ref{eq:BoundOnProbability+1A}) we
obtain the inequality
\begin{equation}
\bP{ U_{r}(\theta) \leq x } \leq {P \choose x} \left ( \frac{ {x
\choose K } }{ {P \choose K} } \right )^{rK}
\label{eq:BoundOnProbability+1}
\end{equation}
from the fact $|\mathcal{P}_x | = {P \choose x}$. Under the
enforced conditions it is the case that
\[
\frac{ {x \choose K} }{ {P \choose K} } = \prod_{\ell=0}^{K-1}
\left ( \frac{x-\ell}{P-\ell} \right ) \leq \left ( \frac{x}{P}
\right )^K
\]
since $\frac{x-\ell}{P-\ell}$ decreases as $\ell$ increases from
$\ell = 0$ to $\ell=K-1$, and the inequality
(\ref{eq:BoundOnProbability+2}) follows by using this fact into
(\ref{eq:BoundOnProbability+1}). \myendpf

The bounds (\ref{eq:BoundOnProbability+2}) trivially hold with
$\bP{ U_{r}(\theta) \leq x } = 0$ when $x =1, \ldots , K-1$ since
we always have $U_r(\theta) \geq K$. We shall make repeated use of
this fact as follows: For all $n,r =1, 2, \ldots , $ with $r < n$,
we have
\begin{eqnarray}
{n \choose r} ~ \bP{ U_{r} (\theta) \leq x } &\leq & {n \choose
r}{P \choose x } \left ( \frac{x}{P} \right)^{rK}
\nonumber \\
&\leq & {\lfloor P/\sigma \rfloor \choose r}{P \choose x } \left (
\frac{x}{P} \right)^{rK} \label{eq:BoundOnProbability+3}
\end{eqnarray}
on the range $x=1, \ldots , \min(rK,P)$ {\em whenever} $\sigma n
\leq P$ for some $\sigma > 0$, a condition needed only for the
last step and which implies $n \leq \lfloor \frac{P}{\sigma}
\rfloor$ since $n$ is an integer.

We are now in a position to improve on the bound
(\ref{eq:CrudeUpperBound}).

\begin{lemma}
{\sl Consider positive integers $K$ and $P$ such that $K \leq P$.
With $n=2,3, \ldots $ and $r=1, \ldots , n$, we have
\begin{eqnarray}
\bP{ A_{n,r} (\theta) } &\leq& \bP{ U_r(\theta) \leq x } e^{-(n-r)
\frac{K^2}{P}}
\nonumber \\
& & + ~ \bP{ C_{r} (\theta) } e^{- (n-r) \frac{K}{P}(x+1) }
\label{eq:BasicDecompositionSumOfPiecesTighter}
\end{eqnarray}
for each positive integer $x$. }
\label{lem:BasicDecompositionSumOfPiecesTighter}
\end{lemma}

\myproof Fix $n=2,3, \ldots $ and pick $r=2, \ldots , n-1$. For
each positive integer $x$, consider the decomposition
\begin{eqnarray}
\bP{ A_{n,r} (\theta) } &=& \bP{ A_{n,r}(\theta) \cap [
U_r(\theta) \leq x ] }
\nonumber \\
& & + ~ \bP{ A_{n,r}(\theta) \cap [ U_r(\theta) > x ] }.
\label{eq:DecompositionForEachR}
\end{eqnarray}

Using (\ref{eq:BoundComputePA_{n,r}}) (with $E=[U_r(\theta) \leq
x]$) and the bound $U_r(\theta) \geq K$, we get
\begin{eqnarray}
\lefteqn{ \bP{ A_{n,r}(\theta) \cap [ U_r(\theta) \leq x ] }
 } & &
\nonumber \\
&\leq& \bP{ C^\star_{r}(\theta) \cap [ U_r(\theta) \leq x ] }
\cdot e^{-(n-r) \frac{K^2}{P} }
\nonumber \\
&\leq& \bP{ U_r(\theta) \leq x } \cdot e^{-(n-r) \frac{K^2}{P} } .
\label{eq:BoundComputePA_{n,r}+1}
\end{eqnarray}
Invoking (\ref{eq:BoundComputePA_{n,r}}) again (this time with
$E=[U_r(\theta) > x]$), we find
\begin{eqnarray}
\lefteqn{ \bP{ A_{n,r}(\theta) \cap [ U_r (\theta) \geq x ] }  } &
&
\nonumber \\
&\leq& \bE{ \1{ C^\star_{r}(\theta) \cap [ U_r(\theta) > x ] }
\cdot e^{-(n-r) \frac{K}{P} \cdot U_{r} (\theta) } }
\nonumber \\
&\leq& \bP{ C_r(\theta) } e^{-(n-r) \frac{K}{P}(x+1) }
\label{eq:BoundComputePA_{n,r}+2}
\end{eqnarray}
since $U_{r} (\theta) \geq x+1 $ on the event $[ U_r(\theta) > x
]$. We complete the proof by combining
(\ref{eq:DecompositionForEachR}),
(\ref{eq:BoundComputePA_{n,r}+1}) and
(\ref{eq:BoundComputePA_{n,r}+2}). \myendpf

\section{Outlining the proof of Proposition \ref{prop:OneLawAfterReduction}}
\label{sec:Outline}

It is now clear how to proceed: Consider a strongly  admissible
scaling $P,K: \mathbb{N}_0 \rightarrow \mathbb{N}_0$ as in the
statement of Proposition \ref{prop:OneLawAfterReduction}. Under
(\ref{eq:AdmissibilityB}) we necessarily have $\lim_{n \rightarrow
\infty } \frac{P_n}{K_n} = \infty $ as discussed at the end of
Section \ref{sec:RoadMap}; see
(\ref{eq:RatioConditionStrong+Consequence1}). As a result,
$\lim_{n \rightarrow \infty} r_n(\theta_n) = \infty$, and for any
given integer $R \geq 2$ we have
\begin{equation}
R < r_n(\theta_n) , \quad n \geq n^\star(R)
\label{eq:SelectingValuesOfNandR}
\end{equation}
for some finite integer $n^\star(R)$.

For the time being, pick an integer $R \geq 2$ (to be specified in
Section \ref{sec:ToShowSecondPiece1}), and on the range $n \geq
n^\star(R)$ consider the decomposition
\begin{eqnarray}
\sum_{r=2}^{\lfloor \frac{n}{2} \rfloor} {n \choose r} ~ \bP{
A_{n,r} (\theta_n) } &=& \sum_{r=2}^{ R } {n \choose r} ~ \bP{
A_{n,r} (\theta_n) }
\label{eq:AnotherDecomposition}\\
& & +~ \sum_{r=R+1}^{ r_n (\theta) } { n \choose r} ~ \bP{ A_{n,r}
(\theta_n) }
\nonumber \\
& & +~ \sum_{ r=r_n (\theta_n)+1}^{\lfloor \frac{n}{2} \rfloor} {
n \choose r} ~ \bP{ A_{n,r} (\theta_n) } . \nonumber
\end{eqnarray}
Let $n$ go to infinity: The desired convergence
(\ref{eq:OneLawToShow}) will be established if we show
\begin{equation}
\lim_{n \rightarrow \infty} \sum_{r=2}^{ R } { n \choose r} ~ \bP{
A_{n,r} (\theta_n) } = 0 , \label{eq:StillToShow0}
\end{equation}
\begin{equation}
\lim_{n \rightarrow \infty} \sum_{r=R+1}^{ r_n (\theta_n) } { n
\choose r} ~ \bP{ A_{n,r} (\theta_n) } = 0 \label{eq:StillToShow1}
\end{equation}
and
\begin{equation}
\lim_{n \rightarrow \infty} \sum_{ r=r_n (\theta_n)+1}^{\lfloor
\frac{n}{2} \rfloor} { n \choose r} ~ ~ \bP{ A_{n,r} (\theta_n) }
= 0 . \label{eq:StillToShow2}
\end{equation}

The next sections are devoted to proving the validity of
(\ref{eq:StillToShow0}), (\ref{eq:StillToShow1}) and
(\ref{eq:StillToShow2}) by repeated applications of Lemma
\ref{lem:BasicDecompositionSumOfPiecesTighter}. We address these
three cases by making use of the bounds
(\ref{eq:BasicDecompositionSumOfPiecesTighter}) with
\[
x = \lfloor (1+\varepsilon) K_n \rfloor, \quad \varepsilon \in
(0,\frac{1}{2} ),
\]
\[
x = \lfloor \lambda r K_n \rfloor, \quad \lambda \in (0,1),
\]
and
\[
x = \lfloor \mu P_n \rfloor, \quad \mu \in (0,1),
\]
respectively. Throughout, we make repeated use of the standard
bounds
\begin{equation}
{n \choose r} \leq \left ( \frac{e n}{r} \right )^r, \quad
\begin{array}{c}
r=1, \ldots , n \\
n=1,2, \ldots \\
\end{array}
\label{eq:CombinatorialBound1}
\end{equation}
Finally, from convexity we note the inequality
\begin{equation}
(x + y )^p \leq 2^{p-1} (x^p + y ^p ), \quad
\begin{array}{c}
x,y \geq 0 \\
p \geq 1.
\end{array}
\label{eq:ConvexityInequality}
\end{equation}

Before getting on the way, we close this section by highlighting
key differences between our approach and the one used in the
papers \cite{BlackburnGerke,
DiPietroManciniMeiPanconesiRadhakrishnan2006}. The observation
yielding (\ref{eq:BasicIdea+UnionBound2}), which forms the basis
of our discussion, is also used in some form as the starting point
in both these references. However, these authors did not take
advantage of the fact that the sufficiently tight bound
(\ref{eq:ProbabilityOfC}) is available for the probability of the
event $C_{r}(\theta)$, a consequence of the {\em exact} expression
(\ref{eq:ProbabilityOfTree}). Through this bound, we can leverage
strong admissibility (via (\ref{eq:EquivalentCondition})) to get
\[
\left ( 1 - q(\theta_n) \right ) \leq (1+\delta) \cdot
\frac{K^2_n}{P_n}
\]
for $n$ sufficiently large with any $0 < \delta < 1$, in which
case
\[
\bP{ C_r(\theta_n) } \leq r^{r-2} \left ( (1+\delta) \cdot
\frac{K^2_n}{P_n} \right )^{r-1}
\]
for each $r=2,3, \ldots , n $. This opens the way to using the
properties of the scaling by means of its deviation function
defined by (\ref{eq:DeviationCondition}) -- Such a line of
arguments cannot be made if the scaling is merely admissible.

The bound (\ref{eq:BasicDecompositionSumOfPiecesTighter}) arises
from the need to efficiently bound the rv $U_r(\theta_n)$. Indeed,
if it were the case that $U_r(\theta_n)= r K_n$ for each $r=1,
\ldots, \lfloor \frac{n}{2} \rfloor$, then the conjecture
(\ref{eq:CONJECTUREScalingK+P})-(\ref{eq:CONJECTUREZeroOneLawConnectivity})
would readily follow as in Erd\H{o}s-R\'enyi graphs by simply
making use of the bound (\ref{eq:CrudeUpperBound}), e.g., see the
arguments in \cite{Bollobas, DraiefMassoulieLMS,
SpencerSaintFlour1991}. In addition, the constraint $U_r(\theta_n)
\leq \min (rK_n, P_n)$ already suggests that the cases $r K_n \leq
P_n$ and $P_n < r K_n$ be considered separately, with a different
decomposition (\ref{eq:BasicDecompositionSumOfPiecesTighter}) on
each range -- This was also the approach taken in the references
\cite{BlackburnGerke,
DiPietroManciniMeiPanconesiRadhakrishnan2006}. Interestingly
enough, a further decomposition of the range $r=1, \ldots, \lfloor
\frac{P_n}{K_n} \rfloor$ is needed to establish Theorem
\ref{thm:ZeroOneLaw+Connectivity}. In particular, using the bound
(\ref{eq:BasicDecompositionSumOfPiecesTighter}) with $x= \lfloor
\lambda r K_n \rfloor$ for sufficiently small $\lambda$ in $(0,1)$
across the entire range $r=1, \ldots, \lfloor \frac{P_n}{K_n}
\rfloor$ would not suffice for very small values of $r$: In that
range the obvious bound $U_r(\theta_n) \geq K_n$ might be tighter
than $U_r(\theta_n) \geq \lfloor \lambda r K_n \rfloor$, and
another form of the bound
(\ref{eq:BasicDecompositionSumOfPiecesTighter}) is needed to
obtain the desired results, hence (\ref{eq:AnotherDecomposition}).

\section{Establishing (\ref{eq:StillToShow0})}
\label{sec:ToShowSecondPiece0}

Consider a strongly admissible scaling $P,K: \mathbb{N}_0
\rightarrow \mathbb{N}_0$ whose deviation function $\alpha :
\mathbb{N}_0 \rightarrow \mathbb{R}$ satisfies $\lim_{n \to
\infty}{\alpha_n}=\infty$. According to this scaling, for each $r
=2,3, \ldots $ and $n=r+1,r+2, \ldots $, replace $\theta$ by
$\theta_n$ in Lemma \ref{lem:BasicDecompositionSumOfPiecesTighter}
with $x= \lfloor(1+\varepsilon) K_n \rfloor$ for some
$\varepsilon$ in $(0, \frac{1}{2})$. For an arbitrary integer
$R\geq 2$, the convergence (\ref{eq:StillToShow0}) will follow if
we show that
\begin{equation}
\lim_{n\rightarrow \infty} {n \choose r} ~ \bP{ C_{r}(\theta_n ) }
e^{- (n-r) \frac{K_n}{P_n} (\lfloor (1+\varepsilon) K_n \rfloor +
1) } = 0 \label{eq:fixed_r_prop_1}
\end{equation}
and
\begin{equation}
\lim_{n\rightarrow \infty} {n \choose r} ~ \bP{ U_{r} (\theta_n)
\leq \lfloor (1+\varepsilon)K_n \rfloor } e^{-(n-r)
\frac{K_n^2}{P_n}} = 0 \label{eq:fixed_r_prop_2}
\end{equation}
for {\em each} $r=2,3, \ldots$. These two convergence statements
are established below in Proposition \ref{prop:fixed_r_prop_1} and
Proposition \ref{prop:fixed_r_prop_2}, respectively.

\begin{proposition}
{\sl Consider a strongly  admissible scaling $P,K: \mathbb{N}_0
\rightarrow \mathbb{N}_0$ whose deviation function $\alpha :
\mathbb{N}_0 \rightarrow \mathbb{R}$ satisfies $\lim_{n \to
\infty}{\alpha_n}=\infty$. With $\varepsilon > 0$, the convergence
(\ref{eq:fixed_r_prop_1}) holds for each $r=2,3, \ldots $. }
\label{prop:fixed_r_prop_1}
\end{proposition}

\myproof Pick $r=2,3, \ldots $ and $\varepsilon > 0$, and consider
a strongly  admissible scaling $P,K: \mathbb{N}_0 \rightarrow
\mathbb{N}_0$. We combine the bounds (\ref{eq:ProbabilityOfC}) and
(\ref{eq:CombinatorialBound1}) to write
\begin{eqnarray}
\lefteqn{ {n \choose r} \bP{ C_r (\theta_n) }
e^{-(n-r)\frac{K_n}{P_n} \left( \lfloor (1+\varepsilon) K_n
\rfloor+1 \right)} } &&
\nonumber \\
&\leq& \left ( \frac{e n}{r} \right )^r r^{r-2} \left (1 -
q(\theta_n) \right )^{r-1} e^{-(n-r)\frac{K_n}{P_n}
\left(\lfloor(1+\varepsilon) K_n\rfloor+1\right)}
\nonumber \\
&\leq& \left ( \frac{e^r}{r^2} \right ) n^r \left (1 - q(\theta_n)
\right )^{r-1} e^{-(n-r) \frac{K_n^2}{P_n} (1+\varepsilon) }
\label{eq:ForEachFixedS+1}
\end{eqnarray}
for all $n=r+1, r+2, \ldots $. Thus, it follows from Lemma
\ref{lem:AsymptoticEquivalence} (via
(\ref{eq:AsymptoticsEquivalence})) that the convergence
(\ref{eq:fixed_r_prop_1}) will be established if we show that
\begin{equation}
\lim_{n \rightarrow \infty} n^r \left ( \frac{K^2_n}{P_n}
\right)^{r-1} e^{-(n-r) \frac{K^2_n}{P_n}(1+\varepsilon) } = 0.
\label{eq:ForEachFixedREquivalent}
\end{equation}
This step relies on the strong admissibility of the scaling.

On the range where (\ref{eq:ForEachFixedS+1}) holds, we find with
the help of (\ref{eq:DeviationCondition}) that
\begin{eqnarray}
\lefteqn{ n^r \left ( \frac{K^2_n}{P_n} \right )^{r-1} e^{-(n-r)
\frac{K^2_n}{P_n}(1+\varepsilon) }} &&
\nonumber \\
&=& n^r  \left (  \frac{\log n+\alpha_n }{n} \right)^{r-1} \cdot
e^{-(n-r) \frac{ \log n+\alpha_n }{n}(1+\varepsilon) }
\nonumber \\
&=& n \left ( \log n + \alpha_n \right )^{r-1} \cdot e^{ - (1 +
\varepsilon)(1 - \frac{r}{n}) \log n } \cdot e^{ - (1 +
\varepsilon)(1 - \frac{r}{n}) \alpha_n }
\nonumber\\
&=& n^{1 - (1 + \varepsilon)(1 - \frac{r}{n})} \cdot \left ( \log
n + \alpha_n \right )^{r-1} \cdot e^{ - (1 + \varepsilon)(1 -
\frac{r}{n}) \alpha_n }
\nonumber \\
&=& n^{- \varepsilon + (1 + \varepsilon)\frac{r}{n}} \cdot \left (
\log n + \alpha_n \right )^{r-1} \cdot e^{ - (1 + \varepsilon)(1 -
\frac{r}{n}) \alpha_n } . \label{eq:ForEachFixedRBound}
\end{eqnarray}

Under the condition $\lim_{n \rightarrow \infty} \alpha_n =
\infty$ it is plain that
\[
\lim_{n \rightarrow \infty} n^{- \varepsilon + (1 +
\varepsilon)\frac{r}{n}} (\log n)^{r-1} e^{ - (1 + \varepsilon)(1
- \frac{r}{n}) \alpha_n } = 0
\]
and
\[
\lim_{n \rightarrow \infty} n^{- \varepsilon + (1 +
\varepsilon)\frac{r}{n}} \alpha_n^{r-1} e^{ - (1 + \varepsilon)(1
- \frac{r}{n}) \alpha_n } = 0.
\]
Letting $n$ go to infinity in (\ref{eq:ForEachFixedRBound}) we
readily get (\ref{eq:ForEachFixedREquivalent}) by making use of
(\ref{eq:ConvexityInequality}). \myendpf

\begin{proposition}
{\sl Consider a strongly admissible scaling $P,K: \mathbb{N}_0
\rightarrow \mathbb{N}_0$ whose deviation function $\alpha :
\mathbb{N}_0 \rightarrow \mathbb{R}$ satisfies $\lim_{n \to
\infty}{\alpha_n}=\infty$. For every $\varepsilon$ in $(0,
\frac{1}{2})$, the convergence (\ref{eq:fixed_r_prop_2}) holds for
each $r=2,3, \ldots $. } \label{prop:fixed_r_prop_2}
\end{proposition}

\myproof Pick $r=2,3, \ldots $ and $\varepsilon $ in $(0,
\frac{1}{2})$, and consider a strongly  admissible scaling $P,K:
\mathbb{N}_0 \rightarrow \mathbb{N}_0$. For $n$ sufficiently
large, we use (\ref{eq:BoundOnProbability+2}) with $x=\lfloor
(1+\varepsilon) K_n \rfloor$ to obtain
\begin{eqnarray}
\lefteqn{ {n \choose r} \bP{ U_r (\theta_n ) \leq \lfloor
(1+\varepsilon) K_n \rfloor } } &&
\nonumber\\
&\leq& {n \choose r} {P_n \choose{\lfloor
K_n(1+\varepsilon)\rfloor}} \left(\frac{\lfloor
K_n(1+\varepsilon)\rfloor}{P_n}\right)^{rK_n}
\nonumber\\
&\leq& n^{r} \left(\frac{e P_n }{\lfloor
K_n(1+\varepsilon)\rfloor}\right)^{\lfloor
K_n(1+\varepsilon)\rfloor} \left(\frac{\lfloor
K_n(1+\varepsilon)\rfloor}{P_n}\right)^{rK_n}
\nonumber\\
&\leq& n^{r}\left(e^{\frac{\lfloor K_n(1+\varepsilon)\rfloor}
                          {rK_n-\lfloor K_n(1+\varepsilon)\rfloor}}
\frac{\lfloor K_n(1+\varepsilon)\rfloor}{P_n}
\right)^{rK_n-\lfloor K_n(1+\varepsilon)\rfloor} . \nonumber
\end{eqnarray}

The condition $r \geq 2$ implies the inequalities
\[
\frac{\lfloor K_n(1+\varepsilon)\rfloor} {rK_n-\lfloor
K_n(1+\varepsilon)\rfloor} \leq
\frac{1+\varepsilon}{r-(1+\varepsilon)} \leq
\frac{1+\varepsilon}{1-\varepsilon}
\]
and
\[
rK_n-\lfloor K_n(1+\varepsilon)\rfloor \geq K_n \left ( r-
(1+\varepsilon) \right ) > 0.
\]
Thus, upon setting
\[
\Gamma(\varepsilon)
:=(1+\varepsilon)e^{\frac{1+\varepsilon}{1-\varepsilon}},
\]
we conclude by strong admissibility (in view of
(\ref{eq:RatioConditionStrong+Consequence1}))  that
$\Gamma(\varepsilon) \cdot \frac{K_n}{P_n}<1$ for all $n$
sufficiently large, whence
\[
e^{\frac{\lfloor K_n(1+\varepsilon)\rfloor}
        {rK_n-\lfloor K_n(1+\varepsilon)\rfloor}}
\frac{\lfloor K_n(1+\varepsilon)\rfloor}{P_n} \leq \Gamma
(\varepsilon) \cdot  \frac{K_n}{P_n} < 1
\]
on that range.

There we can write
\begin{eqnarray}
\lefteqn{ {n \choose r} \bP{ U_r (\theta_n ) \leq  \lfloor
(1+\varepsilon) K_n \rfloor } } &&
\nonumber \\
&\leq& n^{r}\left ( \Gamma(\varepsilon) \cdot \frac{K_n}{P_n}
\right ) ^{rK_n-\lfloor K_n(1+\varepsilon)\rfloor}
\nonumber \\
&\leq& n^{r}\left( \Gamma (\varepsilon) \cdot \frac{K_n}{P_n}
\right)^{K_n(r-1-\varepsilon)}
\nonumber \\
&\leq& n^{r}\left(\Gamma(\varepsilon) \cdot
\frac{K_n}{P_n}\right)^{2(r-1-\varepsilon)}
\label{eq:auxiliary1} \\
&\leq& n^{r} \left( \Gamma(\varepsilon) \cdot \frac{K_n^2}{P_n}
\right)^{2(r-1-\varepsilon)}
\nonumber \\
&=& n^{r} \left( \Gamma(\varepsilon) \cdot \frac{\log n +
\alpha_n}{n}\right)^{2(r-1-\varepsilon)}
\nonumber \\
&=& n^{-r+2+2\varepsilon} \left(\Gamma(\varepsilon) \cdot (\log n
+ \alpha_n)\right)^{2(r-1-\varepsilon) } \label{eq:FirstFactor}
\end{eqnarray}
where we obtain (\ref{eq:auxiliary1}) upon using the fact $K_n
\geq 2$. On the other hand we also have
\begin{eqnarray}
e^{-(n-r)\frac{K^2_n}{P_n} } = e^{-(n-r)\frac{\log n +
\alpha_n}{n} } = n^{-(1 - \frac{r}{n}) } \cdot e^{ - \frac{n-r}{n}
\alpha_n }. \label{eq:SecondFactor}
\end{eqnarray}

Therefore, upon multiplying (\ref{eq:FirstFactor}) and
(\ref{eq:SecondFactor}) we see that Proposition
\ref{prop:fixed_r_prop_1} will follow if we show that
\begin{equation}
\lim_{n \rightarrow \infty} n^{-r+1+2\varepsilon + \frac{r}{n} }
\cdot \left(\log n + \alpha_n\right)^{2(r-1-\varepsilon)} \cdot
e^{ - \frac{n-r}{n} \alpha_n } =0 . \label{eq:Desired1}
\end{equation}
The choice of $\varepsilon$ and $r$ ensures that $r-1-\varepsilon
> 0$ and $-r+1+2\varepsilon + \frac{r}{n} < 0$ for all $n$
sufficiently large. The condition $\lim_{n \rightarrow \infty }
\alpha_n = \infty$ now yields
\begin{equation}
\lim_{n \rightarrow \infty} n^{-r+1+2\varepsilon + \frac{r}{n} }
\cdot \left(\log n \right)^{2(r-1-\varepsilon)} \cdot e^{ -
\frac{n-r}{n} \alpha_n } =0 \label{eq:Desired1+A}
\end{equation}
and
\begin{equation}
\lim_{n \rightarrow \infty} n^{-r+1+2\varepsilon + \frac{r}{n} }
\cdot \alpha_n ^{2(r-1-\varepsilon)} \cdot e^{ - \frac{n-r}{n}
\alpha_n } = 0 . \label{eq:Desired1+B} \end{equation} The desired
conclusion (\ref{eq:Desired1}) follows by making use of
(\ref{eq:Desired1+A}) and (\ref{eq:Desired1+B}) with the help of
the inequality (\ref{eq:ConvexityInequality}). \myendpf

\section{Establishing (\ref{eq:StillToShow1})}
\label{sec:ToShowSecondPiece1}

In order to establish (\ref{eq:StillToShow1}) we will need two
technical facts which are presented in Proposition \ref{prop:A}
and Proposition \ref{prop:B}.

\begin{proposition}
{\sl Consider a strongly admissible scaling $P,K: \mathbb{N}_0
\rightarrow \mathbb{N}_0$ whose deviation function $\alpha :
\mathbb{N}_0 \rightarrow \mathbb{R}$ satisfies $\lim_{n \to
\infty}{\alpha_n}=\infty$. With $0 < \lambda < 1$ and integer $R
\geq 2$, we then have \begin{equation} \lim_{n\rightarrow \infty}
\sum_{r=R+1}^{ \lfloor \frac{n}{2} \rfloor } {n \choose r} ~ \bP{
C_{r}(\theta_n ) } e^{- (n-r) \frac{K_n}{P_n} \left(\lfloor
\lambda r K_n\rfloor+1\right)} = 0 \label{eq:A}
\end{equation}
whenever $\lambda$ and $R$ are selected so that
\begin{equation}
2 < \lambda ( R+1 ). \label{eq:A+Condition}
\end{equation}
} \label{prop:A}
\end{proposition}

Proposition \ref{prop:A} is proved in Section
\ref{sec:ProofPropositionRef{prop:A}}. Next, set
\begin{equation}
C(\lambda;\sigma ) := \left( \frac{e^2}{\sigma} \right)^{\frac{
\lambda }{ 1 - 2 \lambda } }, \quad
\begin{array}{c}
\sigma > 0 \\
0 < \lambda < \frac{1}{2}. \\
\end{array}
\label{eq:C(Lambda)}
\end{equation}

\begin{proposition}
{\sl Consider a strongly admissible scaling $P,K: \mathbb{N}_0
\rightarrow \mathbb{N}_0$ whose deviation function $\alpha :
\mathbb{N}_0 \rightarrow \mathbb{R}$ satisfies $\lim_{n \to
\infty}{\alpha_n}=\infty$. If there exists some $\sigma>0$ such
that (\ref{eq:OneLaw+ConnectivityExtraCondition}) holds for all
$n=1, 2, \ldots$ sufficiently large, then
\begin{equation}
\lim_{n\rightarrow \infty} \sum_{r=1}^{ r_n (\theta_n) } {n
\choose r} ~ \bP{ U_r(\theta_n) \leq  \lfloor \lambda r K_n
\rfloor } e^{-(n-r) \frac{K_n^2}{P_n}} = 0 \label{eq:B}
\end{equation}
whenever $\lambda$ in $(0, \frac{1}{2} )$ is selected small enough
so that
\begin{equation}
\max \left (2\lambda\sigma, \lambda^{1-2\lambda}, \lambda
C(\lambda;\sigma) \right ) < 1 . \label{eq:CONDITIONonLambda}
\end{equation}
} \label{prop:B}
\end{proposition}

A proof of Proposition \ref{prop:B} can be found in Section
\ref{sec:ProofPropositionRef{prop:B}}. Note that for any $\sigma
>0$, $\lim_{\lambda \downarrow 0} \lambda C(\lambda;\sigma) = 0$
and $\lim_{\lambda \downarrow 0} \lambda^{1-2\lambda}=0$, hence
the condition (\ref{eq:CONDITIONonLambda}) can always be met by
suitably selecting $\lambda > 0$ small enough.

We now turn to the proof of (\ref{eq:StillToShow1}): Keeping in
mind Proposition \ref{prop:A} and Proposition \ref{prop:B}, we
select $\lambda$ sufficiently small in $(0, \frac{1}{2} )$ to meet
the condition (\ref{eq:CONDITIONonLambda}) and then pick any
integer $R \geq 2$ sufficiently large to ensure
(\ref{eq:A+Condition}). Next consider a strongly admissible
scaling $P,K: \mathbb{N}_0 \rightarrow \mathbb{N}_0$ whose
deviation function $\alpha : \mathbb{N}_0 \rightarrow \mathbb{R}$
satisfies the condition $\lim_{n \to \infty}{\alpha_n}=\infty$.
Then, for each $n \geq n^\star (R)$ (with $n^\star (R)$ as
specified at (\ref{eq:SelectingValuesOfNandR})), replace $\theta$
by $\theta_n$ according to this scaling, and for each $r = R+1,
\ldots , r_n(\theta_n)$, set $x=\lfloor \lambda r K_n \rfloor$ in
Lemma \ref{lem:BasicDecompositionSumOfPiecesTighter} with
$\lambda$ as specified earlier.

With these preliminaries in place, we see from Lemma
\ref{lem:BasicDecompositionSumOfPiecesTighter} that
(\ref{eq:StillToShow1}) holds if both limits
\[
\lim_{n\rightarrow \infty} \sum_{r=R+1}^{ r_n(\theta_n) } {n
\choose r} ~ \bP{ C_{r}(\theta_n ) } e^{ - (n-r) \frac{K_n}{P_n}
(\lfloor \lambda r K_n \rfloor + 1 ) } = 0
\]
and
\[
\lim_{n\rightarrow \infty} \sum_{r=R+1}^{ r_n (\theta_n) } {n
\choose r} ~ \bP{ U_{r}(\theta_n) \leq \lfloor \lambda r K_n
\rfloor } e^{-(n-r) \frac{K_n^2}{P_n}} = 0
\]
hold. However, under (\ref{eq:A+Condition}) and
(\ref{eq:CONDITIONonLambda}), these two convergence statements are
immediate from Proposition \ref{prop:A} and Proposition
\ref{prop:B}, respectively. \myendpf

\section{Establishing (\ref{eq:StillToShow2})}
\label{sec:ToShowSecondPiece2}

The following two results are needed to establish
(\ref{eq:StillToShow2}). The first of these results is given next
with a proof available in Section
\ref{sec:ProofPropositionRef{prop:D}}.

\begin{proposition}
{\sl Consider a strongly admissible scaling $P,K: \mathbb{N}_0
\rightarrow \mathbb{N}_0$ whose deviation function $\alpha :
\mathbb{N}_0 \rightarrow \mathbb{R}$ satisfies $\lim_{n \to
\infty}{\alpha_n}=\infty$. If there exists some $\sigma > 0$ such
that (\ref{eq:OneLaw+ConnectivityExtraCondition}) holds for all
$n=1,2, \ldots $ sufficiently large, then
\begin{eqnarray}
\lim_{n \rightarrow \infty} \sum_{r= r_n(\theta_n)+1}^{\lfloor
\frac{n}{2} \rfloor} {n \choose r} ~ \bP{ U_r(\theta_n) \leq
\lfloor \mu P_n \rfloor } e^{-(n-r) \frac{K_n^2}{P_n} } = 0
\nonumber
\end{eqnarray}
whenever $\mu$ in $(0,\frac{1}{2})$ is selected so that
\begin{equation}
\max \left ( 2 \left ( \sqrt{\mu} \left ( \frac{e}{ \mu } \right
)^{\mu} \right )^\sigma, \sqrt{\mu} \left ( \frac{e}{ \mu }
\right)^{\mu} \right ) < 1 . \label{eq:ConditionOnMU+1}
\end{equation}
} \label{prop:D}
\end{proposition}

We have $\lim_{\mu \downarrow 0} \left ( \frac{e}{ \mu }
\right)^{\mu} =1$, whence $\lim_{\mu \downarrow 0} \sqrt{\mu}
\left ( \frac{e}{ \mu } \right)^{\mu} = 0$, and
(\ref{eq:ConditionOnMU+1}) can be made to hold for any $\sigma>0$
by taking $\mu > 0$ sufficiently small. The second proposition is
established in Section \ref{sec:ProofPropositionRef{prop:C}}.

\begin{proposition}
{\sl Consider an admissible scaling $P,K: \mathbb{N}_0 \rightarrow
\mathbb{N}_0$ whose deviation function $\alpha : \mathbb{N}_0
\rightarrow \mathbb{R}$ satisfies $\lim_{n \to
\infty}{\alpha_n}=\infty$. If there exists some $\sigma > 0$ such
that (\ref{eq:OneLaw+ConnectivityExtraCondition}) holds for all
$n=1,2, \ldots $ sufficiently large, then
\[
\lim_{n\rightarrow \infty} \sum_{r=r_n(\theta_{n})+1}^{ \lfloor
\frac{n}{2} \rfloor } {n \choose r} ~ \bP{ C_{r}(\theta_n ) } e^{-
(n-r) \frac{K_n}{P_n}(\lfloor \mu P_n \rfloor +1)} = 0
\]
for each $\mu$ in $(0,1)$. } \label{prop:C}
\end{proposition}

The proof of (\ref{eq:StillToShow2}) is now within easy reach:
Consider a strongly admissible scaling $P,K: \mathbb{N}_0
\rightarrow \mathbb{N}_0$ whose deviation function $\alpha :
\mathbb{N}_0 \rightarrow \mathbb{R}$ satisfies $\lim_{n \to
\infty}{\alpha_n}=\infty$. On the range where
(\ref{eq:OneLaw+ConnectivityExtraCondition}) holds, for each $n
\geq n^\star (R)$ (with $n^\star (R)$ as specified at
(\ref{eq:SelectingValuesOfNandR}) where $R$ and $\lambda$ still
satisfy (\ref{eq:A+Condition}) and (\ref{eq:CONDITIONonLambda})),
replace $\theta$ by $\theta_n$ according to this scaling, and set
$x=\lfloor \mu P_n \rfloor$ in Lemma
\ref{lem:BasicDecompositionSumOfPiecesTighter} with $\mu$ as
specified by (\ref{eq:ConditionOnMU+1}). We get
(\ref{eq:StillToShow2}) as a direct consequence of Proposition
\ref{prop:D} and Proposition \ref{prop:C}. \myendpf

\section{A proof of Proposition \ref{prop:A}}
\label{sec:ProofPropositionRef{prop:A}}

Let $\lambda$ and $R$ be as in the statement of Proposition
\ref{prop:A}, and pick a positive integer $n$ such that $2(R + 1)
< n$. Arguments similar to the ones leading to
(\ref{eq:ForEachFixedS+1}) yield
\begin{eqnarray}
\lefteqn{ {n \choose r} ~ \bP{ C_{r}(\theta_n ) } e^{- (n-r)
\frac{K_n}{P_n} \left(\lfloor \lambda r K_n\rfloor+1\right)} } & &
\nonumber \\
&\leq& \left ( \frac{ e^r }{r^2} \right ) n^r e^{-\lambda r (n-r)
\frac{K^2_n}{P_n}} \left (1 - q(\theta_n) \right )^{r-1} \nonumber
\end{eqnarray}
for all $r=1, \ldots , n$. Thus, in order to establish
(\ref{eq:A}), we need only show
\[
\lim_{n\rightarrow \infty} \sum_{r=R+1}^{ \lfloor \frac{n}{2}
\rfloor } \frac{ e^r }{r^2}  n^r e^{-\lambda r (n-r)
\frac{K^2_n}{P_n}} \left (1 - q(\theta_n) \right )^{r-1} = 0 .
\]
As in the proof of Proposition \ref{prop:fixed_r_prop_2}, by the
strong admissibility of the scaling (with the help of
(\ref{eq:AsymptoticEquivalenceWithDelta})), it suffices to show
\begin{equation}
\lim_{n\rightarrow \infty} \sum_{r=R+1}^{ \lfloor \frac{n}{2}
\rfloor } \frac{ e^r }{r^2}  n^r e^{-\lambda r (n-r)
\frac{K^2_n}{P_n}} \left ((1+\delta) \frac{K^2_n}{P_n}
\right)^{r-1} = 0 \label{eq:Intermediary1}
\end{equation}
with $0 < \delta < 1$.

Fix $n=2,3, \ldots$. For each $r=1, \ldots , \lfloor \frac{n}{2}
\rfloor$, we get
\begin{eqnarray}
\lefteqn{ \left ( \frac{ e^r }{r^2} \right ) n^r e^{-\lambda r
(n-r) \frac{K^2_n}{P_n}} \left ((1+\delta) \frac{K^2_n}{P_n}
\right)^{r-1} } &&
\nonumber \\
&=& \left ( \frac{ e^r }{r^2} \right ) n^r e^{-\lambda r (n-r)
\frac{\log n+\alpha_n}{n}} \left ((1+\delta) \frac{\log
n+\alpha_n}{n} \right )^{r-1}
\nonumber \\
&=& n \left ( \frac{ e^r }{r^2} \right ) e^{ -\lambda r (n-r)
\frac{\log n+\alpha_n}{n} } \left ( (1+\delta)( \log n + \alpha_n)
\right )^{r-1}
\nonumber \\
&\leq& n e^r e^{ -\lambda r (1-\frac{r}{n}) \left ( \log
n+\alpha_n \right ) } \left ( (1+\delta)( \log n + \alpha_n)
\right )^{r-1}
\nonumber \\
&\leq& n e^r e^{ -\frac{\lambda}{2} r ( \log n + \alpha_n) } \left
( (1+\delta)( \log n + \alpha_n) \right )^{r-1}
\nonumber \\
&=& n \left ( e^{ 1 -\frac{\lambda}{2} ( \log n + \alpha_n) }
\right )^r \left ( (1+\delta)( \log n + \alpha_n) \right )^{r-1}
\nonumber
\end{eqnarray}
as we note that
\begin{equation}
1-\frac{r}{n} \geq \frac{1}{2}, \quad r=1, \ldots , \left \lfloor
\frac{n}{2} \right \rfloor . \label{eq:BoundOnFraction}
\end{equation}

Next, we set
\[
\Gamma_n (\lambda) :=  n e^{1 -\frac {\lambda}{2}(\log n +
\alpha_n)}
\]
and
\[
a_n (\lambda) :=    e^{1 -\frac{\lambda}{2} ( \log n + \alpha_n )
} (1+\delta) ( \log n +\alpha_n ) .
\]
With this notation we conclude that
\begin{eqnarray}
& & \sum_{r=R+1}^{ \lfloor \frac{n}{2} \rfloor } \left ( \frac{
e^r }{r^2} \right ) n^r e^{-\lambda r (n-r) \frac{K^2_n}{P_n}}
\left ( (1+\delta)\frac{K^2_n}{P_n} \right )^{r-1}
\nonumber \\
&\leq& \Gamma_n (\lambda) \sum_{r=R+1}^{ \lfloor \frac{n}{2}
\rfloor } a_n(\lambda)^{r-1}
\nonumber \\
&\leq& \Gamma_n (\lambda) \sum_{r=R}^{ \infty } a_n(\lambda) ^r.
\label{eq:Intermediary2}
\end{eqnarray}

Obviously, $ \lim_{n \rightarrow \infty} a_n (\lambda)  = 0$ under
the condition $\lim_{n \rightarrow \infty} \alpha_n = \infty$, so
that $a_n (\lambda) < 1 $ for all $n$ sufficiently large. On that
range, the geometric series at (\ref{eq:Intermediary2}) converges
to a finite limit with
\[
\sum_{r=R}^{ \infty } a_n (\lambda)^r = \frac{a_n(\lambda) ^R} {1
- a_n(\lambda) } .
\]
Thus,
\begin{eqnarray}
\lefteqn{ \sum_{r=R+1}^{ \lfloor \frac{n}{2} \rfloor } \left (
\frac{ e^r }{r^2} \right ) n^r e^{-\lambda r (n-r)
\frac{K^2_n}{P_n}} \left ((1+\delta) \frac{K^2_n}{P_n} \right
)^{r-1} } &&
\nonumber \\
&\leq& \Gamma_n (\lambda) \cdot \frac{a_n(\lambda)^R} {1 -
a_n(\lambda) }
\nonumber \\
&=& C_{n,R} (\delta) \cdot n^{1 - \frac{\lambda}{2} (R+1)} \cdot
e^{- \frac{\lambda}{2} (R+1) \alpha_n } \cdot \left ( \log n +
\alpha_n \right )^R \nonumber
\end{eqnarray}
with
\[
C_{n,R} (\delta) := \frac{ e^{R+1} (1+\delta)^R }{1-a_n(\lambda)}.
\]

Under (\ref{eq:A+Condition}), the condition $\lim_{n \rightarrow
\infty} \alpha_n = \infty$ implies
\[
\lim_{n \rightarrow \infty } { n^{1-\frac{\lambda}{2}(R+1)} \cdot
e^{-\frac{\lambda}{2} (R+1) \alpha_n} \cdot \left ( \log n \right
)^R} = 0
\]
and
\[
\lim_{n \rightarrow \infty } n^{1-\frac{\lambda (R+1)}{2}} \cdot
e^{-\frac{\lambda (R+1)}{2}\alpha_n} \cdot \alpha_n ^R = 0.
\]
The desired conclusion (\ref{eq:Intermediary1}) is now immediate
with the help of the inequality (\ref{eq:ConvexityInequality}).
\myendpf

\section{A proof of Proposition \ref{prop:B}}
\label{sec:ProofPropositionRef{prop:B}}

We begin by providing bounds on the probabilities of interest
entering (\ref{eq:B}). Recall the definitions of the quantities
introduced before the statement of Proposition \ref{prop:B}.

\begin{proposition}
{\sl Consider positive integers $K$, $P$ and $n$ such that $2 \leq
K \leq P$ and $\sigma n \leq P$ for some $\sigma>0$. For any
$\lambda $ in $(0, \frac{1}{2})$ small enough to ensure
\begin{equation}
\max \left ( 2 \lambda \sigma , \lambda C(\lambda;\sigma) \right )
< 1, \label{eq:ConditionOnLambda}
\end{equation}
we have
\[
{n \choose r} ~ \bP{ U_{r}(\theta) \leq \lfloor \lambda rK \rfloor
} \leq B(\lambda;\sigma;K )^r \label{eq:BOUND+1}
\]
for all $r=1, \ldots , r_n(\theta)$ where we have set
\[
B(\lambda; \sigma; K) := \max \left ( \lambda^{1-2\lambda},
\lambda^{1-2\lambda} \left ( \frac{e^2}{\sigma} \right
)^{\lambda}, \frac{e^2}{ \sigma K^{K-2}} \right ) .
\label{eq:B(Lambda)}
\]
} \label{prop:BOUND+1}
\end{proposition}

\myproof Pick positive integers $K$, $P$ and $n$ as in the
statement of Proposition \ref{prop:BOUND+1}. For each $r=1,2,
\ldots , n$, we use (\ref{eq:BoundOnProbability+3}) with $x =
\lfloor \lambda r K \rfloor$ to find
\begin{eqnarray}
{n \choose r} ~ \bP{ U_r (\theta) \leq \lfloor \lambda r K \rfloor
} \leq { \lfloor \frac{P}{\sigma} \rfloor \choose r} {P \choose
\lfloor \lambda r K \rfloor } \left ( \frac{ \lfloor \lambda r K
\rfloor }{P} \right )^{rK}. \nonumber
\end{eqnarray}

On the range
\begin{equation}
r= 1, \ldots , r_n(\theta), \label{eq:Range}
\end{equation}
the inequalities
\begin{equation}
r \leq \left \lfloor \frac{P}{K} \right \rfloor - 1 < \frac{P}{K}
\label{eq:Consequences1}
\end{equation}
hold, whence $r < \frac{P}{2}$ since $K \geq 2$. Now if $\lambda$
is selected in $(0 , \frac{1}{2} )$ sufficiently small such that
$2 \lambda \sigma < 1 $, it then follows from
(\ref{eq:Consequences1}) that $\lambda r K < \lambda P <
\frac{P}{2 \sigma}$ so that
\begin{equation}
\lfloor \lambda r K \rfloor \leq \left \lfloor \frac{P}{2\sigma}
\right \rfloor \leq \frac{1}{2} \left \lfloor \frac{P}{\sigma}
\right \rfloor . \label{eq:lambda_r_K_leq_P/2sigma}
\end{equation}
Under these circumstances, we also have
\begin{equation}
rK - \lfloor 2 \lambda r K \rfloor \geq (1 - 2 \lambda ) r K > 0.
\label{eq:PositiveDifference}
\end{equation}
Two possibilities arise:

\noindent {\bf Case I: $r \leq \lfloor \lambda r K \rfloor$ --}
Since $r \leq \lfloor \lambda r K \rfloor \leq \frac{1}{2} \left
\lfloor \frac{P}{\sigma} \right \rfloor $ by
(\ref{eq:lambda_r_K_leq_P/2sigma}), we get
\begin{eqnarray}
\lefteqn{ {n \choose r } ~ \bP{ U_{r} (\theta) \leq \lfloor
\lambda r K \rfloor }} &&
\nonumber \\
&\leq & {\lfloor \frac{P}{\sigma} \rfloor
         \choose \lfloor \lambda r K \rfloor }
        {P \choose \lfloor \lambda r K \rfloor }
\left ( \frac{ \lfloor \lambda r K \rfloor }{P} \right )^{rK}
\nonumber \\
&\leq & \left ( \frac{e \lfloor \frac{P}{\sigma} \rfloor}
                     { \lfloor \lambda rK \rfloor } \right )
^{\lfloor \lambda r K \rfloor} \left ( \frac{eP}{ \lfloor \lambda
rK \rfloor } \right ) ^{\lfloor \lambda r K \rfloor} \left (
\frac{ \lfloor \lambda r K \rfloor }{P} \right)^{rK}
\nonumber \\
&\leq & \left ( \frac{e}{\sigma}\frac{ P}{\lfloor \lambda rK
\rfloor } \right ) ^{\lfloor \lambda r K \rfloor} \left (
\frac{eP}{\lfloor \lambda rK \rfloor } \right ) ^{\lfloor \lambda
r K \rfloor} \left ( \frac{ \lfloor \lambda r K \rfloor }{P}
\right)^{rK}
\nonumber \\
&=& \left(\frac{e^2}{\sigma}\right)^{\lfloor \lambda r K \rfloor }
\left (\frac{ \lfloor \lambda r K \rfloor }{P} \right ) ^{rK - 2
\lfloor \lambda r K \rfloor}
\nonumber \\
&=& \left ( \left(\frac{e^2}{\sigma}\right) ^{\frac{\lfloor
\lambda r K \rfloor }{rK - 2 \lfloor \lambda r K \rfloor } } \cdot
\frac{ \lfloor \lambda r K \rfloor }{P} \right ) ^{rK - 2 \lfloor
\lambda r K \rfloor}
\nonumber \\
&\leq & \left ( \max \left ( 1, C(\lambda;\sigma) \right ) \cdot
\frac{ \lfloor \lambda r K \rfloor }{P} \right )^{rK - 2 \lfloor
\lambda r K \rfloor} \label{eq:CaseI_Inter}
\end{eqnarray}
with $C(\lambda;\sigma )$ given by (\ref{eq:C(Lambda)}) -- In the
last step we made use of (\ref{eq:PositiveDifference}) together
with the fact that
\[
\frac{  \lfloor \lambda r K \rfloor }
        { rK - 2 \lfloor \lambda r K \rfloor }
\leq \frac{  \lambda r K } { rK - 2 \lambda r K } =
\frac{\lambda}{1 - 2 \lambda}
\]
since $\lfloor \lambda r K \rfloor \leq \lambda r K $.

On the range (\ref{eq:Range}), we have $rK \leq P$ from
(\ref{eq:Consequences1}) and substituting this fact into
(\ref{eq:CaseI_Inter}) yields
\[
{n \choose r} ~ \bP{ U_r (\theta) \leq \lfloor \lambda r K \rfloor
} \leq \left ( \lambda \max \left ( 1,  C(\lambda;\sigma) \right )
\right )^{rK - 2 \lfloor \lambda r K \rfloor}.
\]
If $\lambda$ in $(0,\frac{1}{2})$ were selected such that $\lambda
C(\lambda;\sigma) < 1$, then $\lambda \max \left ( 1 ,
C(\lambda;\sigma) \right ) < 1$, and we get
\[
{n \choose r} ~ \bP{ U_r (\theta ) \leq \lfloor \lambda r K
\rfloor } \leq \left ( \lambda \max \left ( 1 , C(\lambda;\sigma)
\right ) \right )^{(1-2\lambda)rK}
\]
by recalling (\ref{eq:PositiveDifference}). With this selection
this last upper bound is largest when $K=1$, whence
\begin{eqnarray}
\lefteqn{ {n \choose r} ~ \bP{ U_r (\theta) \leq \lfloor \lambda r
K \rfloor } } & &
\nonumber \\
&\leq& \left ( \max \left ( \lambda ^{1-2\lambda}, \lambda
^{1-2\lambda} \left( \frac{e^2}{\sigma} \right )^{\lambda} \right)
\right)^r \label{eq:Bound1_for_B_sigma_K}.
\end{eqnarray}

\noindent {\bf Case II: $\lfloor \lambda r K \rfloor \leq r $ --}
On the range (\ref{eq:Range}), we have $\lfloor \lambda r K
\rfloor \leq r \leq \frac{P}{2}$ by virtue of
(\ref{eq:Consequences1}).
This time we find
\begin{eqnarray}
\lefteqn{ {n \choose r } ~ \bP{ U_r (\theta) \leq \lfloor \lambda
r K \rfloor } } & &
\nonumber \\
&\leq & {\lfloor \frac{P}{\sigma} \rfloor \choose r } {P \choose r
} \left ( \frac{ \lfloor \lambda r K \rfloor }{P} \right )^{rK}
\nonumber \\
&\leq & \left ( \frac{e}{r} \left \lfloor \frac{P}{\sigma} \right
\rfloor \right )^r \left ( \frac{eP}{r} \right )^r \left ( \frac{
\lfloor \lambda r K \rfloor }{P} \right )^{rK}
\nonumber \\
&\leq & \left ( \frac{eP}{r \sigma} \right )^r \left (
\frac{eP}{r} \right )^r \left ( \frac{ \lfloor \lambda r K \rfloor
}{P} \right )^{rK} . \nonumber
\end{eqnarray}

The condition $\lfloor \lambda r K \rfloor \leq r $ now implies
\begin{eqnarray}
{n \choose r } ~ \bP{ U_r (\theta) \leq \lfloor \lambda r K
\rfloor } &\leq& \left ( \frac{eP}{r \sigma} \right )^r \left (
\frac{eP}{r} \right )^r \left ( \frac{r}{P} \right )^{rK} .
\nonumber \\
&=& \left ( \frac{e^2}{\sigma} \cdot \left ( \frac{ r }{P}
\right)^{(K- 2)} \right )^r
\nonumber \\
&\leq & \left ( \frac{e^2 }{ \sigma K^{K-2} } \right )^r
\label{eq:Bound2_for_B_sigma_K}
\end{eqnarray}
since $r \leq \frac{P}{K}$ upon using (\ref{eq:Consequences1}).
The proof of Proposition \ref{prop:BOUND+1} is completed by
combining the inequalities (\ref{eq:Bound1_for_B_sigma_K}) and
(\ref{eq:Bound2_for_B_sigma_K}). \myendpf

We can now turn to the proof of Proposition \ref{prop:B}: Consider
positive integers $K$, $P$ and $n$ as in the statement of
Proposition \ref{prop:BOUND+1}. Pick $\lambda $ in $(0,
\frac{1}{2})$ which satisfies (\ref{eq:CONDITIONonLambda}) and
note that (\ref{eq:ConditionOnLambda}) is also valid under this
selection. In the usual manner we get
\begin{eqnarray}
\lefteqn{ \sum_{r=1}^{ r_n (\theta) } {n \choose r} ~ \bP{
U_{r}(\theta) \leq \lfloor \lambda r K \rfloor } \cdot e^{-(n-r)
\frac{K^2}{P}} } &&
\nonumber \\
&\leq& \sum_{r=1}^{ r_n (\theta) } {n \choose r} ~ \bP{
U_r(\theta) \leq \lfloor \lambda r K \rfloor } \cdot e^{-\left(n-
\lfloor \frac{n}{2} \rfloor \right) \frac{K^2}{P}}
\nonumber \\
&\leq& e^{-\frac{n}{2} \frac{K^2}{P}} \sum_{r=1}^{ r_n (\theta) }
{n \choose r} ~ \bP{ U_{r}(\theta) \leq \lfloor \lambda r K
\rfloor }
\nonumber \\
&\leq& e^{-\frac{n}{2} \frac{K^2}{P}} \sum_{r=1}^{ r_n (\theta) }
B(\lambda;\sigma;K)^r \nonumber
\end{eqnarray}
as we invoke Proposition \ref{prop:BOUND+1}. If it is the case
that $B(\lambda;\sigma;K ) <1$, the geometric series is summable
with
\[
\sum_{r=1}^{ r_n (\theta) } B(\lambda;\sigma;K)^r \leq
\sum_{r=1}^{ \infty } B(\lambda;\sigma;K)^r =
\frac{B(\lambda;\sigma;K)}{1-B(\lambda;\sigma;K)},
\]
so that
\begin{eqnarray}
\lefteqn{ \sum_{r=1}^{ r_n (\theta) } {n \choose r} ~ \bP{
U_r(\theta) \leq \lfloor \lambda r K \rfloor } \cdot e^{-(n-r)
\frac{K^2}{P}} } & &
\nonumber \\
&\leq& e^{-\frac{n}{2} \frac{K^2}{P}} \frac{ B(\lambda;\sigma;K)}{
1 - B(\lambda;\sigma;K) }. \label{eq:auxiliary2}
\end{eqnarray}

Now, consider a strongly admissible scaling $P,K: \mathbb{N}_0
\rightarrow \mathbb{N}_0$ whose deviation function $\alpha :
\mathbb{N}_0 \rightarrow \mathbb{R}$ satisfies $\lim_{n \to
\infty}{\alpha_n}=\infty$. On the range where
(\ref{eq:OneLaw+ConnectivityExtraCondition}) holds, replace
$\theta$ by $\theta_n$ in the last inequality according to this
admissible scaling. From (\ref{eq:DeviationCondition}) we see that
\[
K_n^2 = \frac{P_n}{n}(\log n + \alpha_n) \geq \sigma (\log n+
\alpha_n)
\]
so that $\lim_{n \to \infty}{K_n}=\infty$, whence
\[
\lim_{n \to \infty} 
\left (
 \frac{ e^2}{ \sigma K_n^{K_n-2} }
 \right )
=0 .
\]
Moreover, any $\lambda$ in the interval $(0, \frac{1}{2})$
satisfying (\ref{eq:CONDITIONonLambda}) also satisfies the
condition $\lambda C(\lambda; \sigma) < 1 $, so that
\[
\lambda^{1-2\lambda} \left(\frac{e^2}{\sigma}\right)^{\lambda} =
(\lambda C(\lambda; \sigma))^{1-2\lambda} < 1 .
\]

As a result, under (\ref{eq:CONDITIONonLambda}) we see that
\[
\lim_{n \rightarrow \infty} B ( \lambda; \sigma; K_n) = \max \left
( \lambda^{1-2\lambda}, \lambda^{1-2\lambda} \left (
\frac{e^2}{\sigma} \right)^{\lambda} \right ) < 1
\]
whence $ B ( \lambda; \sigma; K_n) < 1 $ for all $n$ sufficiently
large. Therefore, on that range (\ref{eq:auxiliary2}) is valid
under the enforced assumptions with $\theta$ is replaced by
$\theta_n$, and we obtain
\begin{eqnarray}
\nonumber \lefteqn{ \sum_{r=1}^{ r_n (\theta) } {n \choose r} ~
\bP{ U_r (\theta_n) \leq \lfloor \lambda r K_n \rfloor } \cdot
e^{-(n-r) \frac{K_n^2}{P_n}} } &&
\nonumber \\
&\leq& e^{-\frac{n}{2} \frac{ \log n + \alpha_n }{n}} \cdot \left
( \frac{ B(\lambda;\sigma;K_n)}{1- B(\lambda;\sigma;K_n) } \right
)
\nonumber\\
&=&n^{-\frac{1}{2}}e^{-\frac{\alpha_n}{2}} \cdot \left ( \frac{
B(\lambda;\sigma;K_n)}
             {1- B(\lambda;\sigma;K_n) } \right ) .
\nonumber
\end{eqnarray}
Finally, let $n$ go to infinity in this last expression: The
condition $\lim_{n \rightarrow \infty} \alpha_n = \infty$ implies
$\lim_{n \to \infty } n^{-\frac{1}{2}}e^{-\frac{\alpha_n}{2}}=0$
and this completes the proof. \myendpf

\section{A proof of Proposition \ref{prop:D}}
\label{sec:ProofPropositionRef{prop:D}}

Proposition \ref{prop:D} is an easy consequence of the following
bound.

\begin{proposition}
{\sl Consider positive integers $K$ and $P$ such that $2 \leq K $
and $2K \leq P$. For each $\mu$ in $(0, \frac{1}{2})$, we have
\begin{eqnarray}
\lefteqn{ \sum_{r= r_n(\theta)+1}^{\lfloor \frac{n}{2} \rfloor} {n
\choose r} ~ \bP{ U_r(\theta) \leq \lfloor \mu P \rfloor }
e^{-(n-r) \frac{K^2}{P} } } & &
\nonumber \\
&\leq& \left ( 2 e^{- \frac{K^2}{2P} } \right )^n \left (
\sqrt{\mu} \left ( \frac{e}{ \mu } \right )^{\mu} \right)^P
\label{eq:D+Auxiliary}
\end{eqnarray}
for all $n=2,3,\ldots $. } \label{prop:D+Auxiliary}
\end{proposition}

\myproof Fix $n=2,3, \ldots $. In establishing
(\ref{eq:D+Auxiliary}) we need only consider the case $r_n(\theta)
< \lfloor \frac{n}{2} \rfloor $ (for otherwise
(\ref{eq:D+Auxiliary}) trivially holds), so that $r_n(\theta) =
r(\theta)$ and $r_n(\theta) + 1 = \lfloor \frac{P}{K} \rfloor $.
The range $r_n(\theta)+1 \leq r \leq \lfloor \frac{n}{2} \rfloor$
is then equivalent to
\[
\left \lfloor \frac{P}{K} \right \rfloor \leq r \leq \left \lfloor
\frac{n}{2} \right \rfloor ,
\]
hence
\[
r K \geq \left(\frac{P}{K}-1\right)K\geq \frac{P}{2}
\]
as we make use of the condition $2K \leq P$ in the last step.

With $\mu $ in the interval $(0,\frac{1}{2})$ it follows that
\[
\lfloor \mu P \rfloor \leq \frac{P}{2} \leq \min ( rK, P )
\]
and the bound (\ref{eq:BoundOnProbability+2}) applies with $x =
\lfloor \mu P \rfloor $ for all $r= r(\theta)+1, \ldots , \lfloor
\frac{n}{2} \rfloor$.

With this in mind, recall (\ref{eq:BoundOnFraction}). We then get
\begin{eqnarray}
\lefteqn{ \sum_{r= r_n(\theta)+1}^{\lfloor \frac{n}{2} \rfloor} {n
\choose r} ~ \bP{ U_r(\theta) \leq \lfloor \mu P \rfloor }
e^{-(n-r) \frac{K^2}{P} } } &&
\nonumber \\
&\leq& \sum_{r= r(\theta)+1}^{\lfloor \frac{n}{2} \rfloor} {n
\choose r} {P \choose \lfloor \mu P \rfloor } \left ( \frac{
\lfloor \mu P \rfloor }{P} \right )^{rK} e^{-(n-r) \frac{K^2}{P} }
\nonumber \\
&\leq& e^{-\frac{n}{2} \frac{K^2}{P} } \sum_{r=
r(\theta)+1}^{\lfloor \frac{n}{2} \rfloor} {n \choose r} \left (
\frac{e P}{ \lfloor \mu P \rfloor } \right )^{\lfloor \mu P
\rfloor} \left ( \frac{ \lfloor \mu P \rfloor }{P} \right )^{rK}
\nonumber \\
&\leq& e^{-\frac{n}{2} \frac{K^2}{P} } \sum_{r=
r(\theta)+1}^{\lfloor \frac{n}{2} \rfloor} {n \choose r}
e^{\lfloor \mu P \rfloor} \left ( \frac{\lfloor \mu P \rfloor }{P}
\right )^{rK-\lfloor \mu P \rfloor}
\nonumber \\
&\leq& e^{-\frac{n}{2} \frac{K^2}{P} } \sum_{r=
r(\theta)+1}^{\lfloor \frac{n}{2} \rfloor} {n \choose r}
e^{\lfloor \mu P \rfloor} \mu^{rK-\lfloor \mu P \rfloor}
\label{eq:D+AuxiliaryBound} \\
&\leq& e^{-\frac{n}{2} \frac{K^2}{P} } \left ( \frac{e }{ \mu }
\right )^{\lfloor \mu P \rfloor} \left (
\sum_{r=r(\theta)+1}^{\lfloor \frac{n}{2} \rfloor} {n \choose r}
\right ) \mu^{\frac{P}{2}} \nonumber
\end{eqnarray}
since $ \frac{P}{2} \leq rK $ for all $r = r(\theta) +1, \ldots ,
\lfloor \frac{n}{2} \rfloor $ as pointed out earlier. The passage
to (\ref{eq:D+AuxiliaryBound}) made use of the fact that $r K -
\lfloor \mu P \rfloor \geq 0$. The binomial formula now implies
\begin{equation}
\sum_{r= r(\theta)+1}^{\lfloor \frac{n}{2} \rfloor} {n \choose r}
\leq 2^n, \label{eq:Bin}
\end{equation}
so that
\begin{eqnarray}
\lefteqn{ \sum_{r= r_n(\theta)+1}^{\lfloor \frac{n}{2} \rfloor} {n
\choose r} ~ \bP{ U_r(\theta) \leq  \lfloor \mu P \rfloor }
e^{-(n-r)\frac{K^2}{P} }  } & &
\nonumber \\
&\leq& \left ( 2 e^{- \frac{K^2}{2P} } \right )^n \left (
\frac{e}{ \mu  } \right )^{ \mu P } \mu^{\frac{P}{2}} \nonumber
\end{eqnarray}
and the desired conclusion (\ref{eq:D+Auxiliary}) follows.
\myendpf

Now, if in Proposition \ref{prop:D+Auxiliary}, we assume that
$\sigma n \leq P$ for some $\sigma > 0$, then the inequality
\[
\left ( \sqrt{\mu} \left ( \frac{e}{ \mu } \right )^{\mu} \right
)^P \leq \left ( \sqrt{\mu} \left ( \frac{e}{ \mu } \right )^{\mu}
\right )^{\sigma n}
\]
follows as soon as
\begin{equation}
\sqrt{\mu} \left ( \frac{e}{ \mu } \right )^{\mu} < 1,
\label{eq:ConditionOnMU+2}
\end{equation}
and (\ref{eq:D+Auxiliary}) takes the more compact form
\begin{eqnarray*}
\lefteqn{ \sum_{r= r_n(\theta)+1}^{\lfloor \frac{n}{2} \rfloor} {n
\choose r} ~ \bP{ U_r(\theta) \leq \lfloor \mu P \rfloor }
e^{-(n-r) \frac{K^2}{P} } } & &
\nonumber \\
&\leq& \left ( 2 e^{- \frac{K^2}{2P} } \left ( \sqrt{\mu} \left (
\frac{e}{ \mu } \right )^{\mu} \right )^\sigma \right )^n .
\end{eqnarray*}

To conclude the proof of Proposition \ref{prop:D}, observe that
(\ref{eq:ConditionOnMU+2}) is implied by selecting $\mu$ in
$(0,\frac{1}{2})$ according to (\ref{eq:ConditionOnMU+1}). In that
case, consider a strongly admissible scaling $P,K: \mathbb{N}_0
\rightarrow \mathbb{N}_0$. On the range where
(\ref{eq:OneLaw+ConnectivityExtraCondition}) holds, replace
$\theta$ by $\theta_n$ in the last inequality
according to this scaling. This yields
\begin{eqnarray}
\lefteqn{ \sum_{r= r_n(\theta_n)+1}^{\lfloor \frac{n}{2} \rfloor}
{n \choose r} ~ \bP{ U_r(\theta_n) \leq \lfloor \mu P_n \rfloor }
e^{-(n-r) \frac{K_n^2}{P_n} }} &&
\nonumber \\
&\leq& \left ( 2 e^{- \frac{K_n^2}{2P_n} } \left ( \sqrt{\mu}
\left ( \frac{e}{ \mu } \right )^{\mu} \right )^\sigma \right )^n
\nonumber \\
&\leq& \left ( 2 \left ( \sqrt{\mu} \left ( \frac{e}{ \mu }
\right)^{\mu} \right )^\sigma \right )^n  . \nonumber
\end{eqnarray}
Letting $n$ go to infinity in this last inequality, we readily get
the desired conclusion from (\ref{eq:ConditionOnMU+1}). \myendpf

\section{A proof of Proposition \ref{prop:C}}
\label{sec:ProofPropositionRef{prop:C}}

Consider positive integers $K$ and $P$ such that $2 \leq K \leq
P$, and pick $\mu$ in the interval $(0,1)$. For each $n=2,3,
\ldots $, crude bounding arguments yield
\begin{eqnarray}
\lefteqn{ \sum_{r=r_n(\theta)+1}^{ \lfloor \frac{n}{2} \rfloor }
{n \choose r} ~ \bP{ C_{r}(\theta ) } \cdot e^{ - (n-r)\frac{
K}{P} \left( \lfloor \mu P \rfloor +1\right) }} &&
\nonumber \\
&\leq& \sum_{r=r_n(\theta)+1}^{ \lfloor \frac{n}{2} \rfloor }
{n\choose r} ~ e^{ - (n-r)\frac{ K}{P} (\mu P) }
\nonumber\\
&\leq& \left( \sum_{r=r_n(\theta)+1}^{\lfloor \frac{n}{2} \rfloor}
{n \choose r} \right) e^{- \frac{n}{2} K \mu }
\nonumber\\
&\leq& 2^n e^{-\frac{n}{2} K \mu } \label{eq:C+Auxiliary}
\end{eqnarray}
where we have used (\ref{eq:BoundOnFraction}) and (\ref{eq:Bin}).

To complete the proof of Proposition \ref{prop:C}, consider an
admissible scaling $P,K: \mathbb{N}_0 \rightarrow \mathbb{N}_0$
whose deviation function $\alpha : \mathbb{N}_0 \rightarrow
\mathbb{R}$ satisfies $\lim_{n \to \infty}{\alpha_n}=\infty$.
Replace $\theta$ by $\theta_n$ in (\ref{eq:C+Auxiliary}) according
to this admissible scaling so that
\[
\sum_{r= r_n(\theta_n)+1}^{\lfloor \frac{n}{2} \rfloor} {n \choose
r} ~ \bP{ C_r(\theta_n) } e^{-(n-r) \frac{K_n}{P_n} \lfloor \mu
P_n \rfloor } \leq \left ( 2 e^{-\frac{\mu K_n}{2} } \right )^n .
\]

Let $n$ go to infinity in this last inequality: The condition
(\ref{eq:OneLaw+ConnectivityExtraCondition}) implies
\[
K^2_n =  \frac{ \log n + \alpha_n }{n} \cdot P_n \ \geq \sigma
\left ( \log n + \alpha_n \right )
\]
for $n=1,2, \ldots$ sufficiently large, whence $\lim_{n
\rightarrow \infty} K_n = \infty $ under the assumed condition
$\lim_{n \rightarrow \infty} \alpha_n = \infty$. Consequently,
\[
\lim_{n \rightarrow \infty} \left ( 2 e^{- \frac{ \mu K_n }{2} }
\right ) = 0
\]
and the desired conclusion follows. \myendpf

\section*{Acknowledgment}

The authors thank the anonymous reviewers for their careful
reading of the original manuscript;
their comments helped improve the final version of this paper.

\bibliographystyle{IEEE}


\begin{biographynophoto}{Osman Ya\u{g}an} (S'07) received the B.S.
degree in Electrical and Electronics Engineering from the Middle
East Technical University, Ankara (Turkey) in 2007, and the Ph.D
degree in Electrical and Computer Engineering from the University
of Maryland, College Park, MD in 2011.

He was a visiting Postdoctoral Scholar at Arizona State University 
during Fall 2011.
Since December 2011, he has been a Postdoctoral Fellow with CyLab
at Carnegie Mellon University. His research interests
include security in wireless networks, percolation theory, random
graphs and their applications.
\end{biographynophoto}

\begin{biographynophoto}{Armand M. Makowski} (M'83–-SM'94–-F'06)
received the Licence en Sciences Math\'ematiques from the
Universit\'e Libre de Bruxelles in 1975, the M.S. degree in
Engineering-Systems Science from U.C.L.A. in 1976 and the Ph.D.
degree in Applied Mathematics from the University of Kentucky in
1981. In August 1981, he joined the faculty of the Electrical
Engineering Department at the University of Maryland College Park,
where he is Professor of Electrical and Computer Engineering. He
has held a joint appointment with the Institute for Systems
Research since its establishment in 1985.

Armand Makowski was a C.R.B. Fellow of the Belgian-American
Educational Foundation (BAEF) for the academic year 1975-76; he is
also a 1984 recipient of the NSF Presidential Young Investigator
Award and became an IEEE Fellow in 2006.

His research interests lie in applying advanced methods from the
theory of stochastic processes to the modeling, design and
performance evaluation of engineering systems, with particular
emphasis on communication systems and networks.
\end{biographynophoto}

\end{document}